\documentclass[aos,MSNbibl,nameyear,dvips]{arximspdf}
\usepackage{multirow}
\usepackage{graphicx}

\doi{10.1214/11-AOS873}
\volume{39}
\issue{3}
\pubyear{2011}
\firstpage{1689}
\lastpage{1719}

\makeatletter
\newtheorem{theorem}{Theorem}
\newtheorem{lemma}{Lemma}
\newtheorem{assumption}{Assumption}
\makeatother

\begin{document}
\begin{frontmatter}

\title{Testing whether jumps have finite or~infinite~activity}
\runtitle{Finite or infinite jump activity}

\begin{aug}
\author[A]{\fnms{Yacine} \snm{A\"{i}t-Sahalia}\corref{}\thanksref{a1}\ead[label=e1]{yacine@princeton.edu}}
\and
\author[B]{\fnms{Jean} \snm{Jacod}\ead[label=e2]{jean.jacod@upmc.fr}}
\runauthor{Y. A\"{i}t-Sahalia and J. Jacod}
\affiliation{Princeton University and UPMC (Universit\'{e} Paris-6)}
\address[A]{Department of Economics\\
Princeton University and NBER\\
Princeton, New Jersey 08544-1021\\
USA\\
\printead{e1}} 
\address[B]{Institut de Math\'{e}matiques\\ \quad  de Jussieu\\
CNRS UMR 7586\\
Universit\'{e} P. et M. Curie (Paris-6)\\
75252 Paris C\'{e}dex 05\\
France\\
\printead{e2}}
\end{aug}
\thankstext{a1}{Supported in part by NSF Grants DMS-05-32370 and SES-0850533.}

\received{\smonth{3} \syear{2010}}
\revised{\smonth{10} \syear{2010}}

%
\begin{abstract}
We propose statistical tests to discriminate between the finite and
infinite activity of jumps in a semimartingale
discretely observed at high frequency. The two statistics allow for a
symmetric treatment of the problem: we can either
take the null hypothesis to be finite activity, or infinite activity.
When implemented on high-frequency stock returns,
both tests point toward the presence of infinite-activity jumps in the data.
\end{abstract}

%

\end{frontmatter}

\section{Introduction}\label{secintro}\label{intro}

Traditionally, models with jumps in finance have relied on Poisson processes,
as in \citet{merton76}, \citet{balltorous83} and \citet
{bates91}. These
jump-diffusion models allow for a finite number of jumps in a finite time
interval, with the idea that the Brownian-driven diffusive part of the model
captures normal asset price variations while the Poisson-driven jump
part of
the model captures large market moves in response to unexpected information.
More recently, financial models have been proposed that allow for infinitely
many jumps in finite time intervals, using a variety of specifications, such
as the variance gamma model of \citet{madanseneta90} and
\citet{madancarrchang98}, the hyperbolic model of \citet
{eberleinkeller95}, the
CGMY model of \citet{carrgemanmadanyor02} and the finite moment
log stable
process of \citet{carrwu03b}. These models can capture both small
and frequent
jumps, as well as large and infrequent ones.

In this paper, we develop statistical procedures to discriminate empirically
between the two situations of finite and infinite number of jumps, while
allowing in both cases for the presence of a continuous component in the
model. While many theoretical models make use of one or the other type
of\vadjust{\goodbreak}
jumps, no statistical test has been proposed so far that can identify
the most
likely type in a given data set, as existing tests have focused on the
issue of
testing for the presence of jumps but not distinguishing between different
types of jumps; see \citet{yacjf02}, \citet{carrwu03},
\citet{barndorffnielsenshephard04}, \citet{huangtauchen05},
\citet{abd07},
\citet{jiangoomen08}, \citet{leemykland08}, \citet
{yacjacod09a} and
\citet{leehannig10}, among others.

The setup we consider is one where a univariate process $X$ is observed
on a
fixed time interval $[0,T]$, at discretely and regularly spaced times
$i\Delta_{n}$. In a typical high-frequency financial application, $X$
will be
the log of an asset price, the length of observation $T$ ranges from,
say, one
day to one year, while the sampling interval $\Delta_{n}$ is small, typically
measured in seconds. Assuming that the observed path has jumps, we want to
test whether there are a finite number of jumps or not on that path, two
properties commonly referred to as ``finite
activity'' or ``infinite
activity'' for the jump component of $X$.

Our aim is to provide asymptotic testing procedures, as the time lag
$\Delta_{n}$ between successive observations goes to $0$, allowing to decide
which of these two hypotheses is more likely; that is, we want to
decide in
which of the two complementary subsets
%
\begin{eqnarray}
\label{I1}
\Omega_{T}^{f}&=&\{\omega\dvtx t\mapsto X_{t}(\omega)\mbox{ has finitely
many jumps
in }[0,T]\},\nonumber
\\[-8pt]
\\[-8pt]
\Omega_{T}^{i}&=&\{\omega\dvtx t\mapsto X_{t}(\omega)\mbox{ has infinitely
many jumps
in }[0,T]\}
\nonumber
\end{eqnarray}
of the sample space $\Omega$ we are in.

More specifically, we want to find tests with a prescribed asymptotic
significance level, and with asymptotic power going to $1$, to test the null
hypothesis that $\omega$ is in $\Omega_{T}^{f}$, and also the symmetric null
hypothesis that the observed $\omega$ is in $\Omega_{T}^{i}$. We will need
some assumptions on the process $X$, basically that it is an It\^{o}
semimartingale. However, we wish to keep the solution as nonparametric as
possible, and in particular we do not want to specify the structure of the
volatility or of the jumps.

The simple idea behind the two test statistics we propose is common
with our
earlier work on testing whether jumps are present or not, or whether a
continuous component is present. We compute certain power variations of the
increments, suitably truncated and/or sampled at different frequencies.
Related methodologies are being utilized by other authors. For example,
\citet{tauchentodorov10} use the test statistics of \citet
{yacjacod09a}, study
its logarithm for different values of the power argument and contrast the
behavior of the plot above two and below two in order to identify the presence
of a Brownian component. \citet{contmancini} use threshold or
truncation-based
estimators of the continuous component of the quadratic variation, originally
proposed in \citet{mancini01}, in order to test for the presence
of a
continuous component in the price process. The resulting\vadjust{\goodbreak} test is applicable
when the jump component of the process has finite variation, and a test for
whether the jump component indeed has finite variation is also
proposed. In
terms of the Blumenthal--Getoor index $\beta$, this corresponds to testing
whether $\beta<1$.

We aim here to construct test statistics which are simple to compute
and have
the desirable property of being model-free. In particular, no feature
of the
dynamics of the underlying asset price, which can be quite complex with
potentially jumps of various activity levels, stochastic volatility,
jumps in
volatility, etc., need to be estimated in order to compute either the
statistic or its distribution under the null hypothesis. In fact, implementing
the two tests we propose requires nothing more than the computation of various
truncated power truncations.

We consider two testing problems, one where the null hypothesis is
finite jump
activity and its ``dual'' where the null
hypothesis is infinite jump activity. Under the null hypothesis of
finite-activity jumps, the test statistic we propose is similar to the simpler
statistic $S_{n}$ of \citet{yacjacod09a} which was employed to
test for the
presence of jumps, with an additional truncation step. An appropriately
selected truncation mechanism eliminates finite-activity jumps, so that the
probability limit of the statistic $S_{n}$ post-truncation is the same
in this
paper as that of the simpler statistic in the earlier work, under a purely
continuous model. While the result is indeed in that case that
``the answer is the same,'' this is not
completely obvious a priori and still needs to be established mathematically.
And the commonality is limited to probability limits: the two
statistics have
different asymptotic distributions.

Under the reverse scenario, where the null hypothesis is that jumps are
infinitely active, then the statistic $S_{n}^{\prime}$ we propose for this
purpose is radically new and so is its asymptotic behavior. That second
statistic has no relationship to previous work.

As we will see below, when implemented on high-frequency stock returns, both
tests point toward the presence of infinitely active jumps in the data. That
is, in the test where $\Omega_{T}^{f}$ is the null hypothesis, we
reject the
null; in the test where $\Omega_{T}^{i}$ is the null hypothesis, we
fail to
reject the null. This is in line with the empirical results of a companion
paper, \citet{yacjacod09b}, which contains an extension to It\^{o}
semimartingales of the classical Blumenthal--Getoor index $\beta$ for
L\'{e}vy
processes and estimators for $\beta$; see also \citet{belomestny}
for different
estimators. This parameter $\beta$ takes values between $0$ and $2$ and plays
the role of a ``degree of jump
activity'' for
infinitely active jump processes. Then if the estimator of $\beta$ is
found to
be ``high'' in its range $[0,2]$, with a
confidence interval excluding~$0$, as it is the case in the empirical findings
of \citet{yacjacod09b}, it is a strong evidence against finite activity.
However, finite activity implies $\beta=0$, but the converse fails, so using
estimators of $\beta$ can at the best allow for tests when the null is\vadjust{\goodbreak}
``infinite activity,'' and even for
this it
does not allow for determining the asymptotic level of the test. Thus
in fact
the present paper and the other one are complementary, both aiming to
have a
picture as complete as possible of a continuous-time process which is
discretely observed at increasing frequencies. Finally we can also mention
that here the assumptions are significantly weaker than in \citet
{yacjacod09b},
in the sense that the test proposed here is nonparametric, where the estimator
of $\beta$ proposed there is parametric.

The paper is organized as follows. Section~\ref{secmodel} describes our model
and the statistical problem. Our testing procedure is described in Section
\ref{sectests}, and Sections~\ref{secMC} and~\ref{secdata} are
devoted to a simulation study of the tests and an empirical
implementation on
high-frequency stock returns. Technical results are gathered in the
supplemental article [\citet{yacjacod08csupp}].

\section{The model}\label{secmodel}

The underlying process $X$ which we observe at discrete times is a
one-dimensional process which we specify below. Observe that taking a
one-dimensional process is not a restriction in our context since, if
it were
multidimensional, infinitely many jumps on $[0,T]$ means that at least
one of
its components has infinitely many jumps, so the tests below can be performed
separately on each of the components. In all the paper the terminal
time $T$
is fixed. However, it is convenient, and not a restriction, to assume
that the
process $X$ is defined over the whole half-line.

Our structural assumption is that $X$ is an It\^{o} semimartingale on some
filtered space $(\Omega,\mathcal{F},(\mathcal{F}_{t})_{t\geq0},\mathbb{P})$,
which means that its characteristics $(B,C,\nu)$ are absolutely continuous
with respect to Lebesgue measure. $B$ is the drift, $C$ is the quadratic
variation of the continuous martingale part, and $\nu$ is the
compensator of
the jump measure $\mu$ of $X$. In other words, we have
%
\begin{eqnarray}
B_{t}(\omega)&=&\int_{0}^{t}b_{s}(\omega)\,ds, \nonumber
\\[-8pt]
\\[-8pt]C_{t}(\omega)&=&\int_{0}%
^{t}\sigma_{s}(\omega)^{2}\,ds,\qquad\nu(\omega,dt,dx) = dt F_{t}(\omega,dx).
\label{9}%
\nonumber
\end{eqnarray}
Here $b$ and $\sigma$ are optional process, and $F=F_{t}(\omega,dx)$ is a
transition measure from $\Omega\times\mathbb{R}_{+}$ endowed with the
predictable $\sigma$-field into $\mathbb{R}\setminus\{0\}$. One may then
write $X$ as
%
\begin{eqnarray}\label{1}
X_{t} &
=&X_{0}+\int_{0}^{t}b_{s}\,ds+\int_{0}^{t}\sigma_{s}\,dW_{s}\nonumber
\\[-8pt]
\\[-8pt]
&& {} + \underbrace{\int_{0}^{t}\int
x1_{\{|x|\leq1\}}(\mu-\nu)(ds,dx)} _{\mathrm{small\ jumps}}+  \underbrace{\int_{0}^{t}\int
x1_{\{|x|>1\}}\mu(ds,dx)}_{\mathrm{large\ jumps}},
\nonumber
\end{eqnarray}
where $W$ is a standard Wiener process. It is also possible to write
the last
two terms above as integrals with respect to a Poisson measure and its
compensator, but we will not need this here.\vadjust{\goodbreak}

The cutoff level $1$ used to distinguish small and large jumps is arbitrary;
any fixed jump size $\varepsilon>0$ will do. In terms of the definition
(\ref{1}), changing the cutoff level amounts to an adjustment to the
drift of
the process. Ultimately, the question we are asking about the finite or
infinite degree of activity of jumps is a question about the behavior
of the
compensator $\nu$ near $0.$ There are always a finite number of big
jumps. The
question is whether there are a finite or infinite number of small
jumps. This
is controlled by the behavior of $\nu$ near $0$.

\subsection{The basic assumptions}\label{ssBA}

The assumptions we make depend on the null hypothesis we want to test. We
start with a very mild (local) boundedness assumption. Recall that a process
$a_{t}$ is pre-locally bounded if $|a_{t}|\leq n$ for $t<T_{n}$, for a
sequence $(T_{n})$ of stopping times increasing to $+\infty$.

\begin{assumption}
\label{assA1} The processes $b_{t}$, $\sigma_{t}$ and $F_{t}(\mathbb{R}
) 1_{\{F_{t}(\mathbb{R})<\infty\}}$ and $\int(x^{2}\wedge1)F_{t}(dx)$ are
pre-locally bounded.
\end{assumption}

In some cases we will need something more about the drift $b$ and the
volatility~$\sigma$.

\begin{assumption}
\label{assA2} The drift process $b_{t}$ is c\`{a}dl\`{a}g, and the volatility
process $\sigma_{t}$ is an It\^{o} semimartingale satisfying Assumption
\ref{assA1}.
\end{assumption}

Under this assumption we can write $\sigma_{t}$ as (\ref{1}), with a Wiener
process $W^{\prime}$ which may be correlated with $W$. Another (equivalent)
way of writing this is
%
\begin{equation}
\sigma_{t}=\sigma_{0}+\int_{0}^{t}\widetilde{b}_{s}\,ds+\int_{0}^{t}
\widetilde{\sigma}_{s}\,dW_{s}+N_{t}+\sum_{s\leq t}\Delta\sigma_{s}
1_{\{|\Delta\sigma_{s}|>1\}}, \label{2}%
\end{equation}
where $N$ is a local martingale which is orthogonal to the Brownian
motion $W$
and has jumps bounded by $1$. Saying that $\sigma$ satisfies Assumption
\ref{assA1} implies that the compensator of the process $[N,N]_{t}%
+\sum_{s\leq t}1_{\{|\Delta\sigma_{s}|>1\}}$ has the form $\int_{0}^{t}%
n_{s}\,ds$, and the processes $\widetilde{b}_{t}$ and $n_{t}$ are pre-locally
bounded, and the process $\widetilde{\sigma}_{t}$ is c\`{a}dl\`{a}g.

Next, we need conditions on the L\'{e}vy measures $F_{t}$, which are quite
stronger than what is in Assumption~\ref{assA1}. We state here a relatively
restrictive assumption.

\begin{assumption}
\label{assA0} The L\'{e}vy measure $F_{t}=F_{t}(\omega,dx)$ is of the form
%
\begin{eqnarray}
\label{BA1}    F_{t}(dx) &=& \frac{\gamma^{\prime}_{t} (\log(1/|x|))^{\delta
_{t}}%
}{|x|^{1+\gamma_{t}}} \bigl( a_{t}^{(+)}1_{\{0<x\leq z_{t} ^{(+)}\}}%
+a_{t}^{(-)}1_{\{-z_{t}^{(-)}\leq x<0\}} \bigr)\,dx\nonumber
\\[-8pt]
\\[-8pt]
&&{}+F^{\prime}_{t}(dx),
\nonumber
\end{eqnarray}
where, for some pre-locally bounded process $L_{t}\geq1$:
\begin{longlist}[(iii)]
\item[(i)] $a_{t}^{(+)}$, $a_{t}^{(-)}$, $z_{t}^{(+)}$ and $z_{t}^{(-)}$ are
nonnegative predictable processes satisfying
%
\begin{equation}
\label{BA2}\frac{1}{L_{t}}\leq z_{t}^{(+)}\leq1,\qquad\frac
{1}{L_{t}}\leq
z_{t}^{(-)} \leq1,\qquad A_{t}:=a_{t}^{(+)}+a_{t}^{(-)}\leq L_{t},
\end{equation}

\item[(ii)] $\gamma_{t}$, $\gamma^{\prime}_{t}$ and $\delta_{t}$ are predictable
processes, satisfying for some constant $\delta\geq0$
%
\begin{eqnarray}
\label{BA7}0&\leq&\gamma_{t}<2-1/L_{t},\qquad|\delta_{t}|\,\leq L_{t},\qquad
\gamma_{t}=0  \quad \Rightarrow \quad  \delta_{t}=\delta,\nonumber
\\[-8pt]
\\[-8pt]\gamma^{\prime}_{t}&=&
\cases{\displaystyle
\gamma_{t}, & \quad if $\gamma_{t}>0$,\cr
1, & \quad if $\gamma_{t}=0$,
}
\nonumber
\end{eqnarray}

\item[(iii)] $F^{\prime}_{t}=F^{\prime}_{t}(\omega,dx)$ is a signed measure, whose
absolute value $|F^{\prime}_{t}|$ satisfies, for some increasing continuous
function $\phi\dvtx [0,1]\mapsto[0,1]$ with $\phi(0)=0$ and some constant
$c\in[0,1)$:
%
\begin{eqnarray}
\label{BA3}%
\gamma_{t} = 0  \quad \mbox{or} \quad  A_{t}=0, \qquad  x>0 \quad&\Rightarrow&\quad
|F^{\prime}_{t}|([-x,x])\leq L_{t}\phi(x\wedge1),\nonumber
\\[-8pt]
\\[-8pt]
\gamma_{t} >  0  \quad \mbox{and} \quad  A_{t}>0 \quad&\Rightarrow&\quad \int
(|x|^{c\gamma_{t}}\wedge1) |F^{\prime}_{t}|(dx)\leq L_{t}.
\nonumber
\end{eqnarray}
\end{longlist}
\end{assumption}

Equivalently, one could take $\gamma_{t}^{\prime}=1$ identically,
provided in
(\ref{BA2}) $A_{t}\leq L_{t}$ is substituted with $A_{t}(1_{\{\gamma_{t}
=0\}}+\frac{1}{\gamma_{t}} 1_{\{\gamma_{t}>0\}})\leq L_{t}$.

Since $F^{\prime}_{t}$ is allowed to be a signed measure, (\ref{BA1}) does
\textit{not}
mean that $F_{t}(dx)$ restricted to $(0,z_{t}^{(+)}]$, say, has the density
$f_{t}^{(+)}(x)=a_{t}^{(+)}\gamma^{\prime}_{t} (\log(1/|x|))^{\delta_{t}
}/|x|^{1+\gamma_{t}}$; it simply means that the ``leading part'' of
$F_{t}$ on
a small interval $(0, \varepsilon]$ has a density behaving as $f_{t}^{(+)}(x)$
as $x\downarrow0$, and likewise on the negative side.

In all models with jumps of which we are aware in financial economics,
such as
those cited in the first paragraph of the \hyperref[intro]{Introduction}, the L\'{e}vy measure
has a density around $0$, which behaves like $\alpha_{t}^{(\pm)}%
(\log(1/|x|))^{\delta_{t}}/|x|^{1+\gamma_{t}}$ as $x\downarrow0$ or
$x\uparrow0$ (in most cases with $\gamma_{t}$ and $\delta_{t}$
constant). Thus
all these models satisfy Assumption~\ref{assA0}. For instance, it is
satisfied if the discontinuous part of $X$ is a stable process of index
$\beta\in(0,2)$, with $\gamma_{t}=\beta$ and $\delta_{t}=0$ and
$z_{t}^{(\pm
)}=1$, and $a_{t}^{(\pm)}$ being constants; in this case the residual measure
$F_{t}^{\prime}$ is the restriction of the L\'{e}vy measure to the complement
of $[-1,1]$, and (\ref{BA3}) holds for any $c\in(0,1)$. When the discontinuous
part of $X$ is a tempered stable process the assumption is also
satisfied with
the same processes as above, but now the residual measure $F_{t}^{\prime
}$ is
not positive in general, although it again satisfies (\ref{BA3}) with any
$c\in(0,1)$. Gamma and two-sided Gamma processes also satisfy this assumption,
take $\gamma_{t}=0$ and $\delta_{t}=0$ and $z_{t}^{(\pm)}=1$, and $a_{t}
^{(\pm)}$ being constant.

This assumption also accounts for a stable or tempered stable or Gamma process
with time-varying intensity, when $\gamma_{t}$, $\delta_{t}$ and $z_{t}%
^{(\pm)}$ are as above, but $a_{t}^{(\pm)}$ are genuine processes. It also\vadjust{\goodbreak}
accounts for a stable with time-varying index process, as well as for $X$
being the sum of a stable or tempered stable process with jump activity index
$\beta$ plus another process whose jumps have activity strictly less than
$\beta$. Furthermore, any process of the form $Y_{t} = Y_{0}+\int_{0}%
^{t}w_{s}\,dX_{s}$ satisfies Assumption~\ref{assA0} as soon as $X$ does and
$w_{t}$ is locally bounded and predictable. As is easily checked (see Section~1 of the supplemental article [\citet
{yacjacod08csupp}] for a
formal proof), under Assumption~\ref{assA0} the set $\Omega_{T}^{i}$ of
(\ref{I1}) is (almost surely)
%
\begin{equation}
\Omega_{T}^{i} = \{\overline{A}_{T}>0\}\qquad\mbox{where }  \overline{A}
_{t} = \int_{0}^{t}A_{s}\,ds. \label{I11}%
\end{equation}

The previous assumption is designed for the test for which the null is
``finite activity.'' For the symmetric test,
the assumption we need is stronger:

\begin{assumption}
\label{assA00} We have Assumption~\ref{assA0} with $\gamma_{t}=\beta$ [a
constant in $(0,2)$], and $\delta_{t}=0$.
\end{assumption}

The reason we need a stronger assumption under the null of infinite jump
activity is that the asymptotic distribution of the test statistic
under the
null is now driven by the behavior of $F_{t}$ near $0,$ whereas in the
previous situation where the null has finite jump activity it is the Brownian
motion that becomes the driving process for the behavior of the statistic
under the null.

Assumptions~\ref{assA0} and~\ref{assA00} have the advantage of being easily
interpretable and also easy to check for any concrete model. But as a matter
of fact it is possible to substantially weaken them, and we do this in the
next subsection. The reader who is satisfied with the degree of
generality of
Assumption~\ref{assA0} can skip the next subsection and go directly to the
description of the tests in Section~\ref{sectests}.

\subsection{Some weaker assumptions}\label{ssWA}

For a better understanding of what follows, let us first recall the
notion of
Blumenthal--Getoor index (in short, BG index). There are two distinct notions.
First, we have a (random) global BG index over the interval $[0,t]$ defined
by
%
\begin{equation}
\Gamma_{t} = \inf\biggl(p>0\dvtx \int_{0}^{t}ds\int(|x|^{p}\wedge1)
F_{s}(dx)<\infty\biggr).
\label{BG12}%
\end{equation}
This is a nondecreasing $[0,2]$-valued process. It is not necessarily
right-continuous, nor left-continuous, but it is always optional.
Second, we
have an instantaneous BG index $\gamma_{t}$, which is the BG index of the
L\'{e}vy measure $F_{t}$, defined as the following (random) number:
%
\begin{equation}
\gamma_{t} = \inf\biggl(p>0\dvtx \int(|x|^{p}\wedge1) F_{t}(dx)<\infty\biggr), \label
{BG1}%
\end{equation}
which necessarily belongs to $[0,2]$. As a process, $\gamma_{t}$ is
predictable. The symmetrical tail function\vadjust{\goodbreak} $\overline{F}_{t}(x)=F_{t}%
((-x,x)^{c})$ of $F_{t}$ (defined for $x>0$) satisfies for all $\omega$ and
$t$:
%
\begin{eqnarray}
p>\gamma_{t}  \quad &\Rightarrow& \quad  \lim_{x\rightarrow0}x^{p} \overline{F}%
_{t}(x)=0,\nonumber
\\[-8pt]
\\[-8pt] p<\gamma_{t}  \quad &\Rightarrow& \quad  \limsup_{x\rightarrow0}%
x^{p} \overline{F}_{t}(x)=\infty. \label{BG2}%
\nonumber
\end{eqnarray}
In the latter case, $x^{p}\overline{F}_{t}(x)$ does not necessarily converge
to $\infty$ when $x\rightarrow0$. When $x^{p}\overline
{F}_{t}(x)\rightarrow
\infty$ as $x\rightarrow0$ for all $p<\gamma_{t}$, or if $F_{t}=0$, we say
that the measure $F_{t}$ is \textit{regular}: this is the case when,
for example,
the function $\overline{F}_{t}$ is regularly varying at~$0$.

The connections between these two indices are not completely straightforward;
they are expressed in the next lemma, where $\Delta X_{s}$ denotes the
jump of
$X$ at time~$s$:

\begin{lemma}
\label{L22} Outside a $\mathbb{P}$-null set, we have for all $t$:
%
\begin{equation}
\label{BG10}\Gamma_{t} = \inf\biggl(p>0\dvtx \sum_{s\leq t}|\Delta
X_{s}|^{p}<\infty\biggr).
\end{equation}
Moreover if $\lambda$ denotes the Lebesgue measure we have, outside a
$\mathbb{P}$-null set again,
%
\begin{equation}
\label{BG11}\Gamma_{t}(\omega) \geq \lambda- \operatorname{\operatorname
{ess}\operatorname{sup}}\bigl(\gamma
_{s}(\omega)\dvtx s\in[0,t]\bigr),
\end{equation}
and this inequality is an equality as soon as $\sup_{s\in[0,t],x\in
(0,1]}x^{\gamma_{s} +\varepsilon} \overline{F}_{s}(x)<\infty$ for all
$\varepsilon>0$.
\end{lemma}

Our general assumption involves two functions with the following properties
[$\phi$~is indeed like in (\ref{BA3})]:
%
\begin{equation}\label{BG16}
\left\{
\begin{tabular}{p{317pt}}
 \mbox{$\phi\dvtx  [0,1]\rightarrow\lbrack0,\infty)$
is continuous increasing with
$\phi(0)=0$},\\[5pt]
\mbox{$\psi\dvtx  (0,1]\rightarrow\lbrack1,\infty)$
is continuous decreasing with
$\psi(0):=\lim_{x\to0}\psi(x)$}\\
  \mbox{being either $1$ [then $\psi(x)=1$ for
all $x$], or $+\infty$,  in which case}\\
  \mbox{$x^\varepsilon\psi(x)=0$ as
$x\to0$ for all $\varepsilon>0$.}
\end{tabular}
\right. \hspace*{-33pt}
\end{equation}

\begin{assumption}
\label{assA3} The global BG index $\Gamma_{t}$ takes its values in $[0,2)$,
and there are a constant $a\in(0,1]$, two functions $\phi$ and $\psi$
as in
(\ref{BG16}) and pre-locally bounded processes $L(\varepsilon)$ and
$L^{\prime}(p)$, such that for all $\varepsilon>0$ and $p\geq2$ and
$\mathbb{P}\otimes\lambda$-almost all $(\omega,t)$ we have
%
\begin{eqnarray}
\label{BG4}
 \qquad &\displaystyle x\in(0,1]   \quad \Rightarrow \quad   x^{\Gamma_{t}+\varepsilon}\overline
{F}_{t}(x)\leq
L(\varepsilon)_{t},&\nonumber\\
 \qquad &\displaystyle x,y\in(0,1]  \quad \Rightarrow \quad  \overline{F}_{t}(x)-\overline{F}_{t}\bigl(x(1+y)\bigr)\leq
\cases{\displaystyle
\frac{L(\varepsilon)_{t}}{x^{\Gamma_{t}+\varepsilon}}\bigl(y+x^{a(\Gamma
_{t}+\varepsilon)}\bigr) ,\cr\qquad  \mbox{if   $\Gamma_{t}>0$},\vspace*{3pt}\cr\displaystyle%
\frac{L(\varepsilon)_{t}}{x^{\Gamma_{t}+\varepsilon}}\bigl(y+x^{\Gamma
_{t}+\varepsilon}\phi(x)\bigr) ,\cr\qquad  \mbox{if   $\Gamma_{t}=0$},
}&
\\
 \qquad &\displaystyle 0<u\leq1, \qquad  \Gamma_{t}=0  \quad \Rightarrow \quad  \int_{\{|x|\leq u\}
}|x|^{p}F_{t}(dx)\leq
L^{\prime}(p)_{t}u^{p}\psi(u).&\nonumber
\end{eqnarray}
\end{assumption}

This assumption looks complicated, but it is just a mild local boundedness
assumption, which is made even weaker by the fact that we use the
global BG
index $\Gamma_{t}$ instead of the (perhaps more natural) instantaneous
index~$\gamma_{t}$.\looseness=1

Now we introduce the set $\Omega_{T}^{ii}$ which will be the
alternative for
our first test when the null is ``finitely many jumps.'' This set has a
somewhat complicated description, which goes as follows, and for which we
introduce the notation (for $p\geq2$):
%
\begin{equation}
\label{BA5}G(p,u)_{t} = \frac1{u^{p}\psi(u)} \int_{\{|x|\leq u\}}|x|^{p}
F_{t}(dx).
\end{equation}
Then we set
%
\begin{eqnarray}\label{BG14}
   \Omega_{T}^{ii} = \Omega_{T}^{i,\Gamma>0}\cup\Omega_{T}^{i,\Gamma
=0}\qquad\mbox{where}\\
\eqntext{\Omega_{T}^{i,\Gamma>0}  =\mbox{the set on which }   \Gamma_{T}%
>0  \mbox{ and, for all }  a^{\prime}\in(0,1),}\\[-25pt]\nonumber
\end{eqnarray}
\begin{eqnarray}
\label{BG140}
\lambda \Bigl(\Bigl\{s\in[0,T]\dvtx  \lim_{x\to0}x^{a^{\prime} \Gamma_{T}}
\overline{F}_{s}(x)=\infty\Bigr\} \Bigr) > 0, \hspace*{-131pt}&&
\\
\Omega_{T}^{i,\Gamma=0}  =\mbox{the set on which }   \Gamma_{T}%
=0  \mbox{ and, for all }  p\geq2,\hspace*{-131pt}&&\nonumber\\
\label{BG141}
\lambda \Bigl(\Bigl\{s\in[0,T]\dvtx  \liminf_{u\to0} G(p,u)_{s}%
 > 0 \Bigr\}\Bigr) > 0. \hspace*{-130pt}&&
\end{eqnarray}

We will see in Section 1 of the supplemental article
[\citet{yacjacod08csupp}] that $\Omega_{T}^{ii}\subset\Omega
_{T}^{i}$, but
equality may fail. In view of Lemma~\ref{L22}, and if the measures
$F_{t}$ are
regular [see after (\ref{BG2})], the set $\Omega_{T}^{i,\Gamma>0}$ is
equal to
$\{\Gamma_{T}>0\}$, which is also the set where there are infinitely many
jumps due to a positive BG index. The interpretation of the set $\Omega
_{T}^{i,\Gamma=0}$ is more delicate: observe that $G(p,u)_{t}\leq
L^{\prime
}(p)_{t}$ by (\ref{BG11}); then saying that $\omega\in\Omega
_{T}^{i,\Gamma=0}$
amounts to saying that for ``enough''
values of
$t\leq T$ the variables $G(p,u)_{t}$ are not small.

The following is proved there:

\begin{lemma}
\label{L16} Assumption~\ref{assA0} implies Assumption~\ref{assA3}, and we
have $\Omega_{T}^{ii}=\Omega_{T}^{i}$ a.s.
\end{lemma}

As for Assumption~\ref{assA00}, it can be weakened as follows:

\begin{assumption}
\label{assA4} There are two constants $\beta\in(0,2)$ and $\beta^{\prime
}%
\in[0,\beta)$ and a pre-locally bounded process $L_{t}$, such that the
L\'{e}vy
measure $F_{t}$ is of the form
%
\begin{equation}
\label{5}F_{t}(dx) = \frac{\beta}{|x|^{1+\beta}} \bigl( a_{t}^{(+)}%
1_{\{0<x\leq z_{t} ^{(+)}\}}+a_{t}^{(-)}1_{\{-z_{t}^{(-)}\leq x<0\}
} \bigr)\,dx+F_{t}^{\prime}(dx)
\end{equation}
[the same as (\ref{BA1}) with $\delta_{t}=0$ and $\gamma_{t}=\beta$]  with
(\ref{BA2}) and the (signed) measure $F^{\prime}_{t}$ satisfying
%
\begin{equation}
\int(|x|^{\beta^{\prime}}\wedge1) |F^{\prime}_{t}|(dx)\leq L_{t}. \label
{7}%
\end{equation}
\end{assumption}

We associate with this assumption the following increasing process and~set:
%
\begin{equation}
\overline{A}_{t} = \int_{0}^{t}A_{s}\,ds,\qquad\Omega_{T}^{i\beta}%
 = \{\overline{A}_{T}>0\}, \label{BG6}%
\end{equation}
where the exponent in $\Omega_{T}^{i\beta}$ stands for ``infinite activity for the jumps associated with the part of the L\'{e}vy
measure having index of jump activity~$\beta$,'' and we have
$\Omega_{T}^{i\beta}\subset\Omega_{T}^{i}$, up to a null set. Again,
the set
$\Omega_{T}^{i}\setminus\Omega_{T}^{i\beta}$ is not necessarily empty.

Assumption~\ref{assA00} obviously implies Assumption~\ref{assA4}, with the
same $\beta$ and with $\beta^{\prime}=c\beta$, and in this case $\Omega
_{T}^{i\beta}=\Omega_{T}^{i}$. However, Assumption~\ref{assA4}, which is
exactly the assumption under which the estimation of the BG index is performed in
\citet{yacjacod09b}, does not require the measure $F^{\prime}_{t}$
to be finite
when $A_{t}=0$, so it is weaker than Assumption~\ref{assA00}.

\section{The two tests}\label{sectests}

\subsection{Defining the hypotheses to be tested}

Ideally, we would like to construct tests in the following two situations:
%
\begin{eqnarray}
\label{T-100}
&\displaystyle H_{0}\dvtx \Omega_{T}^{f} \quad \mbox{vs.} \quad H_{1}\dvtx \Omega_{T}^{i},&\nonumber
\\[-8pt]
\\[-8pt]
&\displaystyle H_{0}\dvtx \Omega_{T}^{i}\quad \mbox{vs.} \quad H_{1}\dvtx \Omega_{T}^{f}.&%
\nonumber
\end{eqnarray}
In order to derive the asymptotic distributions of the test statistics which
we construct below, we need to slightly restrict these testing hypotheses.
Besides the sets defined in (\ref{BG14}) and (\ref{BG6}), we also
define two
other complementary sets:
%
\begin{equation}
\Omega_{T}^{W}=\biggl \{ \int_{0}^{T}\sigma_{s}^{2}\,ds>0 \biggr\} ,\qquad
\Omega_{T}^{\mathrm{no}W}= \biggl\{ \int_{0}^{T}\sigma_{s}^{2}\,ds=0
\biggr\} .
\label{6}%
\end{equation}
That is, $\Omega_{T}^{W}$ is the set on which the continuous martingale part
of $X$ is not degenerate over $[0,T]$, and the exponents in $\Omega_{T}^{W}$
and $\Omega_{T}^{noW}$ stand for ``the Wiener process is
present,'' and ``no Wiener process is
present,'' respectively.

We will provide tests for the following assumptions (below we state the
assumptions needed for the null to have a test with a given asymptotic level,
and those needed for the alternative if we want the test to have asymptotic
power equal to $1$):
%
\fontsize{10pt}{\baselineskip}\selectfont
\makeatletter
\def\tagform@#1{\normalsize\maketag@@@{(\ignorespaces#1\unskip\@@italiccorr)}}
\makeatother
\begin{eqnarray}\label{T-101}
&&\displaystyle H_{0}\dvtx   \Omega_{T}^{f}\cap\Omega_{T}^{W} \mbox{ (Assumptions~\ref{assA1} and~\ref{assA2})} \quad \mbox{vs.}\quad H_{1}\dvtx
\cases{
\Omega_{T}^{i} \mbox{ (Assumptions~\ref{assA1} and~\ref{assA0})},\cr\displaystyle
\Omega_{T}^{ii} \mbox{ (Assumptions~\ref{assA1} and~\ref{assA3})},
}
\nonumber\hspace*{-25pt}
\\[-8pt]
\\[-8pt]
&&\displaystyle H_{0}\dvtx
\cases{
\Omega_{T}^{i\beta} \mbox{ (Assumptions~\ref{assA1} and~\ref{assA4})},\cr
\Omega_{T}^{i} \mbox{ (Assumptions~\ref{assA1} and~\ref{assA00})}
}
  \quad \mbox{vs.} \quad   H_{1}\dvtx   \Omega_{T}^{f}\cap\Omega_{T}^{W}
 \mbox{ (Assumptions~\ref{assA1}).}
\nonumber\hspace*{-25pt}
\end{eqnarray}
\normalsize
Note that this lets aside the two sets $\Omega_{T}^{f}\cap\Omega
_{T}^{\mathrm{no}W}$ and $\Omega_{T}^{i}\setminus\Omega_{T}^{i\beta}$
in the
null hypothesis: in the context of high-frequency data, no semimartingale
model where $\Omega_{T}^{f}\cap\Omega_{T}^{\mathrm{no}W}$ is not empty
has been
used, to our knowledge. Indeed, on this set the path of $X$ over the interval
$[0,T]$ is a pure drift plus finitely many jumps. It also lets aside
the set
$\Omega_{T}^{i}\setminus\Omega_{T}^{i\beta}$, under Assumptions~\ref{assA1}
and~\ref{assA4}, which is more annoying. Note that the set $\Omega_{T}^{f}$
contains those $\omega$ for which $X(\omega)$ is continuous on $[0,T]$,
although we would not test against infinite activity if we did not know
beforehand that there were some jumps.

Next, we specify the notion of testing when the null and alternative
hypotheses are families of possible outcomes. Suppose now that we want
to test
the null hypothesis ``we are in a subset $\Omega_{0}%
$'' of $\Omega$, against the alternative
``we
are in a subset $\Omega_{1}$,'' with of course $\Omega
_{0}%
\cap\Omega_{1}=\varnothing$. We then construct a critical (rejection) region
$C_{n}$ at stage $n$, that is, when the time lag between observations is
$\Delta_{n}$. This critical region is itself a subset of $\Omega$, which
should depend only on the observed values of the process $X$ at stage
$n$. We
are not really within the framework of standard statistics, since the two
hypotheses are themselves random.

We then take the following as our definition of the asymptotic size:
%
\begin{equation}
a = \sup \Bigl( \limsup_{n} \mathbb{P}(C_{n}\mid A)\dvtx  A\in\mathcal{F}%
, A\subset\Omega_{0} \Bigr) . \label{T-1}%
\end{equation}
Here $\mathbb{P}(C_{n}\mid A)$ is the usual conditional probability knowing
$A$, with the convention that it vanishes if $\mathbb{P}(A)=0$. If
$\mathbb{P}(\Omega_{0})=0$, then $a=0$, which is a natural convention
since in
this case we want to reject the assumption whatever the outcome $\omega
$ is.
Note that $a$ features a form of uniformity over all subsets $A\subset
\Omega_{0}$. As for the asymptotic power, we define it as
%
\begin{equation}
P = \inf \Bigl( \liminf_{n} \mathbb{P}(C_{n}\mid A)\dvtx  A\in\mathcal{F}%
, A\subset\Omega_{1}, \mathbb{P}(A)>0 \Bigr) . \label{T-2}%
\end{equation}
Again, this is a number.

Clearly, and as in all tests in high-frequency statistics, at any given stage
$n$ it is impossible to distinguish between finitely many or infinitely many
jumps (or, for that matter, between a continuous and a discontinuous
path). So
testing such hypotheses can only have an asymptotic meaning, as the mesh
$\Delta_{n}$ goes to $0$. Now, our definition of the asymptotic level
is the
usual one, apart from the fact that we test a given family of outcomes rather
than a given family of laws. For the asymptotic power, it is far from the
typical usual statistical understanding.\vadjust{\goodbreak} Namely saying that $P=1$, as will
often be the case below, does not mean anything like the infimum of the power
over all possible alternatives is $1$; it is rather a form of
consistency on
the set of alternatives.

\subsection{Truncated power variations}

Before stating the results, we introduce some notation, to be used throughout.
We introduce the observed increments of $X$ as
%
\begin{equation}
\Delta_{i}^{n}X = X_{i\Delta_{n}}-X_{(i-1)\Delta_{n}}, \label{T-3}%
\end{equation}
to be distinguished from the (unobservable) jumps of the process,
$\Delta
X_{s}=X_{s}-X_{s-}$. In a typical application, $X$ is a log-asset
price, so
$\Delta_{i}^{n}X$ is the recorded log-return over $\Delta_{n}$ units of time.

We take a sequence $u_{n}$ of positive numbers, which will serve as our
thresholds or cutoffs for truncating the increments when necessary, and will
go to $0$ as the sampling frequency increase. There will be
restrictions on
the rate of convergence of this sequence, expressed in the form
%
\begin{equation}
\sup_{n} \Delta_{n}^{\rho}/u_{n} < \infty\label{M-7}%
\end{equation}
for some $\rho>0$: this condition becomes weaker when $\rho$ increases. Two
specific values for $\rho$, in connection with the power $p>2$ which is used
below, are of interest to us:
%
\begin{equation}
\rho_{1}(p) = \frac{p-2}{2p},\qquad\rho_{2}(p) = \frac{p-2}{4p-4} \wedge
\frac{2p-4}{11p-10}. \label{M-6}%
\end{equation}
These quantities increase when $p$ increases [on $(2,\infty)$], and
$\rho
_{1}(p)>\rho_{2}(p)>0$.

Finally, with any $p>0$ we associate the increasing processes
%
\begin{equation}
B(p,u_{n},\Delta_{n})_{t}=\sum_{i=1}^{[t/\Delta_{n}]}|\Delta_{i}^{n}%
X|^{p} 1_{\{|\Delta_{i}^{n}X|\leq u_{n}\}} \label{T-4}%
\end{equation}
consisting of the sum of the $p$th absolute power of the increments
of $X,$ truncated at level $u_{n}$, and sampled at time intervals
$\Delta_{n}
$. These truncated power variations, used in various combinations, will
be the
key ingredients in the test statistics we construct below.

\subsection{The finite-activity null hypothesis}

We first set the null hypothesis to be finite activity, that is, $\Omega
_{0}=\Omega_{T}^{f}\cap\Omega_{T}^{W}$, whereas the alternative is
$\Omega
_{1}=\Omega_{T}^{i}$.

We choose an integer $k\geq2$ and a real $p>2$. We then propose the test
statistic, which depends on $k$, $p$, and on the truncation level
$u_{n}$, and
on the time interval $[0,T]$, defined as follows:
%
\begin{equation}
S_{n}=\frac{B(p,u_{n},k\Delta_{n})_{T}}{B(p,u_{n},\Delta_{n})_{T}}.
\label{T-5}%
\end{equation}
That is, we compute the truncated power variations at two different
frequencies in the numerator and denominator, but otherwise use the
same power
$p$ and truncation level\vadjust{\goodbreak} $u_{n}.$ Since $k$ is an integer, both truncated
power variations can be computed from the same data sample. If the original
data consist of log-returns sampled every $\Delta_{n}$ units of time, then
sampling every $k\Delta_{n}$ units of time involves simply retaining every
$k$th observation in that same sample.

The first result gives the limiting behavior of the statistic $S_{n}$, in
terms of convergence in probability:

\begin{theorem}
\label{TT1}
(\textup{a}) Under Assumption~\ref{assA1} and if the sequence $u_{n}$
satisfies (\ref{M-7}) with some $\rho<1/2$, we have
%
\begin{equation}
S_{n} \stackrel{\mathbb{P}}{\longrightarrow} k^{p/2-1}
\qquad\mbox{{on the set }} \Omega_{T}^{f}\cap\Omega_{T}^{W}. \label{T-6}%
\end{equation}

(\textup{b}) Under Assumptions~\ref{assA1} and~\ref{assA0} (resp., Assumption~\ref{assA3})
and if
the sequence $u_{n}$ satisfies (\ref{M-7}) with some $\rho<\rho
_{1}(p)$, we
have
%
\begin{equation}
S_{n} \stackrel{\mathbb{P}}{\longrightarrow} 1
\qquad\mbox{on the set }
\Omega_{T}^{i}\ (\mbox{{resp.,} }  \Omega_{T}^{ii}). \label{T-7}%
\end{equation}
\end{theorem}

As the result shows, the statistic $S_{n}$ behaves differently
depending upon
whether the number of jumps is finite or not. Intuitively, if the
number of
jumps is finite, then at some point along the asymptotics the truncation
eliminates them and the residual behavior of the truncated power
variation is
driven by the continuous part of the semimartingale. More specifically,
$B(p,u_{n},\Delta_{n})_{T}$ is of order $O_{p}(\Delta_{n}^{p/2-1}A(p)_{T})$
where $A(p)_{T} = m_{p}\int_{0}^{T}|\sigma_{s}|^{p}\,ds$ is the continuous
variation of order $p$ and $m_{p}$ is a constant defined below. It follows
from this that the ratio in $S_{n}$ has the limit given in (\ref{T-6}) since
the numerator and denominator tend to zero but at different rates.

By setting $\rho<1/2$ in the truncation rate, we are effectively
retaining all
the increments of the continuous part of the semimartingale, and so we indeed
obtain the ``full'' continuous variation
$A(p)_{T}$ after truncation$.$ Note that by contrast, the untruncated power
variation converges when $p>2$ to the discontinuous variation of order $p$,
namely $B(p)_{T} = \sum_{s<T}|\Delta X_{s}|^{p}$, and so we would not have
been able to distinguish finite or infinite jump activity without truncation,
as long as jumps (of any activity) are present.

If the jumps have infinite activity, on the other hand, that is, under the
alternative hypothesis, the asymptotic behavior of the truncated power
variation is driven by the small jumps, whether the Brownian motion is present
or not, and the truncation rate matters. We will see, for example,
that, under
Assumption~\ref{assA00} or~\ref{assA4}, the truncated power variation
is of
order $O_{p}(u_{n}^{p-\beta}).$ Since here we are truncating at the
same level
$u_{n}$ for the two sampling frequencies $\Delta_{n}$ and $k\Delta
_{n},$ this
fact has no consequence on the behavior of the ratio $S_{n}$, which
tends to
$1$ as stated by (\ref{T-7}) even when Assumption~\ref{assA4} fails.

Theorem~\ref{TT1} implies that for the test at hand an a priori reasonable
critical region is $C_{n}=\{S_{n}<c_{n}\}$, for $c_{n}$ between $1$ and
$k^{p/2-1}$, and in particular if we choose $c_{n}=c$ in the interval
$(1,k^{p/2-1})$ the asymptotic level and power are, respectively, $0$ and
$1$ if
the model satisfies Assumption~\ref{assA3}.\vadjust{\goodbreak}

For a more refined version of this test, with a prescribed level $a\in(0,1)$,
we need a central limit theorem associated with the convergence in
(\ref{T-6}), and\vspace*{-2pt} a standardized version goes as follows [we use
$\stackrel
{\mathcal{L}-(s)}{\longrightarrow}$ to denote the stable convergence in law;
see, e.g., \citet{jacodshiryaev2003} for this notion]:

\begin{theorem}
\label{TTCLT1} Assume Assumptions~\ref{assA1} and~\ref{assA2}, and that the
sequence $u_{n}$ satisfies (\ref{M-7}) with some $\rho<1/2$. Then
%
\begin{equation}\quad
(S_{n}-k^{p/2-1})/\sqrt{V_{n}} \stackrel{\mathcal{L}-(s)}{\longrightarrow
} \mathcal{N}(0,1)\qquad\mbox{in restriction to }\Omega_{T}%
^{f}\cap\Omega_{T}^{W}, \label{eqCLTS36}%
\end{equation}
where
%
\begin{eqnarray}\label{T-8}
 \qquad V_{n}  & =& N(p,k) \frac{B(2p,u_{n},\Delta_{n})_{T}}{(B(p,u_{n},\Delta
_{n})_{T})^{2}},\\
\label{T-9}
N(p,k) & =&\frac{1}{m_{2p}}  \bigl( k^{p-2}(1+k)m_{2p}+k^{p-2}(k-1)m_{p}%
^{2}-2k^{p/2-1}m_{k,p} \bigr) , %
\end{eqnarray}
and $m_{p,k}=\mathbb{E}(|U|^{p}|U+\sqrt{k-1} V|^{p})$ and $m_{p}%
=\mathbb{E}(|U|^{p})$ for $U$ and $V$ two independent $\mathcal{N}(0,1)$ variables.
\end{theorem}

Note that to implement this asymptotic distribution in practice, and
hence the
test, we simply need to compute truncated power variations for various powers
(specifically $p$ and $2p$, in this case) and at two different sampling
frequencies ($\Delta_{n}$ and $k\Delta_{n}$). No other estimation is required.
In particular, we do not need to estimate any aspect of the dynamics of the
$X,$ such as its drift or volatility processes or its jump measure. In that
sense, the test statistic is nonparametric, or model-free.

We are now ready to exhibit a critical region for testing $H_{0}\dvtx \Omega
_{T}^{f}\cap\Omega_{T}^{W} $ vs. $H_{1}\dvtx \Omega_{T}^{i}$, or $\Omega_{T}^{ii}$,
or $\Omega_{T}^{i,\Gamma>0}$ (depending on the assumptions), using $S_{n}$
with a prescribed asymptotic level $a\in(0,1)$. Denoting by $z_{a}$ the
$a$-quantile of $\mathcal{N}(0,1)$, that is, $\mathbb{P}(U>z_{a})=a$ where $U$ is
$\mathcal N(0,1)$, we set
%
\begin{equation}
C_{n} = \bigl\{S_{n}<k^{p/2-1}-z_{a}\sqrt{V_{n}}\bigr\}. \label{T-10}%
\end{equation}

More precisely, we can state the level and power of the test under selected
alternatives as follows:

\begin{theorem}
\label{TT2} Under Assumptions~\ref{assA1} and~\ref{assA2}, and if the
sequence $u_{n}$ satisfies (\ref{M-7}) with some $\rho<1/2$, the asymptotic
level of the critical region defined by (\ref{T-10}) for testing the null
hypothesis ``finite jump activity'' (i.e.,
$\Omega_{T}^{f}\cap\Omega_{T}^{W}$) equals $a$. Moreover the asymptotic power
of this test is $1$ in the following cases:
\begin{itemize}
\item under Assumption~\ref{assA0} with $\delta>0$, for the alternative
$\Omega_{T}^{i}$,

\item under Assumption~\ref{assA3} with $\psi(0)=\infty$, for the alternative
$\Omega_{T}^{ii}$,

\item under Assumption~\ref{assA3}, for the alternative $\Omega_{T}%
^{i,\Gamma>0}$.\vadjust{\goodbreak}
\end{itemize}
\end{theorem}

Of course this test has the same asymptotic level $a$ for testing
$\Omega
_{T}^{f}\cap\Omega_{T}^{W}$ against its complement $(\Omega_{T}^{f}\cap
\Omega_{T}^{W})^{c}$, under Assumptions~\ref{assA1} and~\ref{assA2};
but in
this case the asymptotic power is probably not $1$, and may well be $0$.

\subsection{The infinite-activity null hypothesis}

In the second case, we assume Assumption~\ref{assA4} (or Assumption~\ref{assA00}) and
we set the null hypothesis to be infinite activity, that is, $\Omega_{0}=
\Omega_{T}^{i\beta}$ (or $\Omega_{T}^{i}$), whereas the alternative is
$\Omega_{1}=\Omega_{T}^{f}\cap\Omega_{T}^{W}$. Unfortunately, we cannot simply
use the same test statistic $S_{n}$ we proposed in (\ref{T-5}) for this
second testing problem. The reason is that, while the distribution of $S_{n}$
is model-free under the null hypothesis of the first testing problem,
it is no
longer model-free under the null hypothesis of the second testing
problem: its
distribution when jumps have infinite activity depends upon the degree
of jump
activity, $\beta.$ While it is possible to estimate $\beta$
consistently [see
\citet{yacjacod09b}], it would be preferable to construct a
statistic whose
implementation does not require a preliminary estimate of the degree of
jump activity.

In order to design a test statistic which is model-free under the null of
infinite activity, we choose three reals $\gamma>1$ and $p^{\prime
}>p>2$ and
then define a family of test statistics as follows:
%
\begin{equation}
S_{n}^{\prime} = \frac{B(p^{\prime},\gamma u_{n},\Delta_{n})_{T}%
B(p,u_{n},\Delta_{n})_{T}}{B(p^{\prime},u_{n},\Delta_{n})_{T}
B(p,\gamma
u_{n},\Delta_{n})_{T}}. \label{T-11}%
\end{equation}
In other words, unlike the previous statistic $S_{n}$, we now play with
different powers $p$ and $p^{\prime},$ and different levels of truncation
$u_{n}$ and $\gamma u_{n},$ but otherwise sample at the same frequency
$\Delta_{n}$. Once more, the first result states the limiting behavior of
$S_{n}^{\prime}$ in terms of convergence in probability:

\begin{theorem}
\label{TT3} Assume Assumptions~\ref{assA1} and~\ref{assA00} (resp.,
Assumptions~\ref{assA1} and~\ref{assA4}).

\begin{longlist}[(b)]
\item[(a)] If the sequence $u_{n}$ satisfies (\ref{M-7}) with some $\rho<\rho_{1}(p)$,
we have
%
\begin{equation}
S_{n}^{\prime} \stackrel{\mathbb{P}}{\longrightarrow} \gamma^{p^{\prime}
-p}\qquad\mbox{{on the set }} \Omega_{T}^{i}
\ (\mbox{{resp.,} }  \Omega_{T}^{i\beta}). \label{T-12}%
\end{equation}

\item[(b)] If the sequence $u_{n}$ satisfies (\ref{M-7}) with some $\rho<1/2$, we
have
%
\begin{equation}
S_{n}^{\prime} \stackrel{\mathbb{P}}{\longrightarrow} 1\qquad\mbox{{on
the set }} \Omega_{T}^{f}\cap\Omega^{W}_{T}. \label{T-13}%
\end{equation}
\end{longlist}
\end{theorem}

That is, as was the case in Theorem~\ref{TT1}, the test statistic
$S_{n}^{\prime}$ tends to $1$ under the alternative and to a value different
from $1$ under the null. Intuitively, under the alternative of finite jump
activity, the behavior of each one of the four truncated power
variations in
(\ref{T-11}) is driven by the continuous part of the semimartingale. The
truncation level is such that essentially all the Brownian increments are
kept. Then the truncated power variations all tend to zero at rate
$O_{p}(\Delta_{n}^{p/2-1})$ and $O_{p}(\Delta_{n}^{p^{\prime}/2-1})$,
respectively, and by construction the (random) constants of proportionality
$A(p)_{T}$ and $A(p^{\prime})_{T}$ cancel out in the ratios, producing the
limit $1$ given in (\ref{T-13}).\vadjust{\goodbreak}

If, on the other hand, jumps have infinite activity, then the small
jumps are
the ones that matter and the truncation level becomes material,
producing four
terms that all tend to zero but at the different orders\vspace*{-1pt}
$O_{p}(u_{n}^{p-\beta
}),$ $O_{p}(u_{n}^{p^{\prime}-\beta}),$ $O_{p}((\gamma u_{n})^{p-\beta
})$ and
$O_{p}((\gamma u_{n})^{p^{\prime}-\beta})$, respectively, resulting in the
limit (\ref{T-12}). By design, the limit is independent of $\beta.$

A reasonable critical region is $C_{n}=\{S_{n}^{\prime}<c_{n}\}$ with $c_{n}$
between $1$ and $\gamma^{p^{\prime}-p}$, and if $c_{n}=c$ is in the interval
$(1,\gamma^{p^{\prime}-p})$ the asymptotic level and power are,
respectively,
$0$ and $1$. For a test with a prescribed level we again need a standardized
central limit theorem associated with the convergence in (\ref{T-12}):

\begin{theorem}
\label{TTCLT2} Assume Assumptions~\ref{assA1} and~\ref{assA00} with $c<1/2$
(resp., Assumption~\ref{assA4} with $\beta^{\prime}<\beta/2$), and if the sequence
$u_{n}$ satisfies (\ref{M-7}) with some $\rho<\rho_{2}(p)$, we have
%
\begin{equation}
 \qquad (S_{n}^{\prime}-\gamma^{p^{\prime}-p})/\sqrt{V_{n}^{\prime}} \stackrel
{\mathcal{L}-(s)}{\longrightarrow} \mathcal{N}(0,1) \qquad\mbox{{in
restriction to }} \Omega_{T}^{i} \ (\mbox{{resp.,} }  \Omega
_{T}^{i\beta}),
\label{eqCLTS43}%
\end{equation}
where
%
\begin{eqnarray}
\label{T-14}
V_{n}^{\prime} & =&\gamma^{2p^{\prime}-2p}\biggl (\frac{B(2p,u_{n},\Delta
_{n})_{T}}{(B(p,u_{n},\Delta_{n})_{T})^{2}} +(1-2\gamma^{-p})\frac
{B(2p,\gamma
u_{n},\Delta_{n})_{T}}{(B(p,\gamma u_{n},\Delta_{n})_{T})^{2}}\nonumber\\
&&\hphantom{\gamma^{2p^{\prime}-2p}\biggl (}{} +\frac{B(2p^{\prime},u_{n},\Delta_{n})_{T}}{(B(p^{\prime},u_{n},\Delta
_{n})_{T})^{2}}+(1-2\gamma^{-p^{\prime}})\frac{B(2p^{\prime},\gamma
u_{n},\Delta_{n})_{T}}{(B(p^{\prime},\gamma u_{n},\Delta_{n})_{T})^{2}%
}\nonumber
\\[-8pt]
\\[-8pt]
&&\hphantom{\gamma^{2p^{\prime}-2p}\biggl (}{} -2\frac{B(p+p^{\prime},u_{n},\Delta_{n})_{T}}{B(p,u_{n},\Delta_{n}
)_{T} B(p^{\prime},u_{n},\Delta_{n})_{T}}\nonumber\\
&&\hspace*{5pt}\hphantom{\gamma^{2p^{\prime}-2p}\biggl (}{} -2(1-\gamma^{-p}-\gamma^{-p^{\prime}})\frac{B(p+p^{\prime},\gamma
u_{n},\Delta_{n})_{T}}{B(p,\gamma u_{n},\Delta_{n})_{T} B(p^{\prime},
\gamma
u_{n},\Delta_{n})_{T}} \biggr).
\nonumber
\end{eqnarray}
\end{theorem}

As was the case for $S_{n}$ in Theorem~\ref{TTCLT1}, the asymptotic
distribution of $S_{n}^{\prime}$ under the null is again model-free.
While the
expression (\ref{T-14}) looks complicated, implementing it simply requires
the computation of truncated power variations $B$ of order $p,$
$p^{\prime}$,
$p+p^{\prime},$ $2p$ and $2p^{\prime}$, and at truncation levels
$u_{n}$ and
$\gamma u_{n}.$ No other aspects of the dynamics of the $X$ process,
such as
its degree of jump activity $\beta$, for instance, need to be estimated.

The critical region for testing $H_{0}\dvtx \Omega_{T}^{i}$ or $\Omega
_{T}^{i\beta
} $ vs. $H_{1}\dvtx \Omega_{T}^{f}$ using $S_{n}^{\prime}$ with a prescribed
asymptotic level $a\in(0,1)$ will be
%
\begin{equation}
C_{n}^{\prime} = \bigl\{S_{n}^{\prime}<\gamma^{p^{\prime}-p}-z_{a}\sqrt
{V_{n}^{\prime}}\bigr\}. \label{T-15}%
\end{equation}
We can state more precisely the level and power of the test as follows:

\begin{theorem}
\label{TT4} Under Assumptions~\ref{assA1} and~\ref{assA00} with $c<1/2$
(resp., Assumption~\ref{assA4} with $\beta^{\prime}<\beta/2$), and if the sequence
$u_{n}$ satisfies (\ref{M-7}) with some $\rho<\rho_{2}(p)$, the asymptotic\vadjust{\goodbreak}
level of the critical region defined by (\ref{T-15}) for testing the null
hypothesis ``infinite activity for
jumps''
(i.e., $\Omega_{T}^{i}$, resp., $\Omega_{T}^{i\beta}$, against $\Omega
_{T}^{f} \cap\Omega_{T}^{W}$) equals $a$. Moreover the asymptotic power of
this test is $1$.
\end{theorem}

Under the null hypothesis the rate of convergence is $1/u_{n}^{\beta/2}$
(contrary to the situation of Theorem~\ref{TT2}, where the rate was
$1/\sqrt{\Delta_{n}}$ whatever $\beta$ and $u_{n}$ were). So, although
asymptotically $u_{n}$ does not explicitly show in the test itself, one should
probably take $u_{n}$ as small as possible (we have no choice as to
$\beta$,
of course). That is, one should take $\rho$ as large as possible, which in
turn results in choosing $p$ as large  as possible [recall
(\ref{M-6})]. We discuss actual choices in practice in Sections~\ref{secMC}
and~\ref{secdata} below.

\subsection{Microstructure noise}\label{subsecnoise}

In practice, the tests presented here need a lot of data to be
effective, that
is, we need a high sampling frequency. However, in this case, it is well
established that the so-called \textit{microstructure noise} may be a relevant
consideration, and in some cases may even dominate at ultra high frequencies.
It is outside the scope of this paper to provide a complete analysis of what
happens when noise is present, or to contemplate constructing effective
testing procedures in the presence of noise. However, as a first step,
it may
be enlightening to determine at least the limiting behavior (in probability)
of our test statistics in the presence of noise since this may help
guide the
interpretation of the empirical results when the test is implemented in
practice.

We start first with pure additive noise, which is the type of noise considered
by much of the literature, primarily for reasons of tractability,
although it
may not account very well for the microstructure noise encountered in
practice. It gives, however, an insight on what can happen in the
presence of
noise. In this situation, at any given observation time $t$ we actually
observe log-returns perturbed by a noise term $\varepsilon_{t}$ in
%
\begin{equation}
Z_{t} = X_{t}+\varepsilon_{t}. \label{N0}%
\end{equation}
To avoid intricate statements, we make a few basic, and mild,
assumptions on
the noise process $\varepsilon_{t}$: those variables are i.i.d., centered,
independent of the underlying process, with the following property:
%
\begin{equation}
 \label{N2}%
\begin{tabular}{@{}p{298pt}@{}}
 if $s\neq t$ the variable $\varepsilon_t-\varepsilon_s$ has a density
$f$, continuous and positive at $0$.
\end{tabular}
\end{equation}
Note that no moment condition is required, except $\mathbb
{E}(\varepsilon
_{t})=0$, because we consider only truncated increments.

We write $\overline{B}(p,u_{n},\Delta_{n})$ for the variables
introduced in
(\ref{T-4}) if we replace $X_{i\Delta_{n}}$ by $Z_{i\Delta_{n}}$, and likewise
for $\overline{S}_{n}$ and $\overline{S}^{\prime}_{n}$ if we do the same
substitution in (\ref{T-5}) and (\ref{T-11}) (those are associated with some
integer $k$ and some numbers $p,p^{\prime}>2$ and $\gamma>1$). Under
the above
assumption on the noise, we have the following:\vadjust{\goodbreak}

\begin{theorem}
\label{T1} Let the sequence $u_{n}$ satisfy (\ref{M-7}) with some $\rho
<\frac13$. Then, under Assumption~\ref{assA1} we have
%
\begin{equation}
\label{N3}\overline{S}_{n} \stackrel{\mathbb{P}}{\longrightarrow}%
 \frac1k,\qquad\overline{S}_{n}^{\prime} \stackrel{\mathbb
{P}}{\longrightarrow
} \gamma^{p^{\prime}-p}.
\end{equation}
\end{theorem}

It is remarkable that the assumptions in this theorem are much weaker
than in
the previous results, as far as $X$ is concerned. This is of course because
the noise, when present, becomes the prevalent factor. This has the important
consequence that the statistics $\overline{S}_{n}$ and $\overline{S}%
_{n}^{\prime}$ are no longer able to distinguish between the two
hypotheses of
finite or infinite activity when noise dominates. When we use $\overline
{S}_{n}$ we get a limit which differs from both limits in Theorem~\ref{TT1},
so when we test for the null hypothesis of finite activity and the empirical
value of $\overline{S}_{n}$ is close to $1/k$ we can in principle tell that
this is due to the noise, whereas if it is close to $k^{p/2-1}$ or to
$1$ this
is probably due to finite activity, or infinite activity. When we use
$\overline{S}_{n}^{\prime}$ the situation is worse, because noise plus finite
activity leads to accept the hypothesis that we have infinite activity whether
the jumps have finite or infinite activity.

Alternatively, it might be closer to the reality to model the microstructure
noise as a pure rounding noise. There, instead of observing $X_{i\Delta_{n}},$
we observe $[X_{i\Delta_{n}}]_{\alpha},$ that is, $X_{i\Delta_{n}}$
rounded to
the nearest multiple of $\alpha$, where $\alpha$ is the tick size: typically
$\alpha=1$ cent for a decimalized stock, or $\alpha=1/32{\mathrm{nd}}$
of a
dollar, for bond prices. Let us denote by $\widehat{B}(p,u_{n},\Delta_{n})_{t}$
the variables of (\ref{T-4}) when we replace $X_{i\Delta_{n}}$ by
$[X_{i\Delta_{n}}]_{\alpha}$. Theorems~\ref{TT1} and~\ref{TT3} were
based on
limit theorems for the truncated power variations $B(p,u_{n},\Delta_{n})_{t}$.
In the rounded case we indeed have something radically different: as
soon as
$u_{n}<\alpha$, we eliminate \textit{all} increments, because
increments are
multiples of $\alpha$. Therefore $\widehat{B}(p,u_{n},\Delta
_{n})_{t}=0$ for
all $n$ large enough, and of course the ratios (\ref{T-5}) and (\ref{T-11})
make no sense. It follows that in the case of a rounding noise, the statistics
proposed in this paper are totally meaningless. And as a matter of
fact, even
when the whole (rounded) path of $X$ is observed, we cannot decide
whether we
have finitely many or infinitely many jumps. This means that not only our
statistics do not work for this problem, but there cannot exist asymptotically
consistent tests for this problem.

Naturally, these two idealized descriptions of microstructure noise do not
exhaust the possibilities for modeling the noise. One can, for example,
use a
mixed model which mixes additive noise and rounding, or more general forms,
as in \citet{jacodlimyklandpodolskijvetter09}, for example. At
present, however,
it is not clear how our statistics theoretically behave in these more general
cases, nor how to construct asymptotically consistent tests, nor even
if such
tests exist at all.

\section{Simulation results}\label{secMC}%

We now report simulation results documenting the finite sample
performance of
the test statistics $S_{n}$ and $S_{n}^{\prime}$ in finite samples
under their
respective null and alternative hypotheses.\vadjust{\goodbreak} We calibrate the values to be
realistic for a liquid stock trading on the NYSE. We use observation lengths
ranging from $T=1$ day, consisting of $6.5$ hours of trading, that is,
$n=23\mbox{,}400$ seconds, to $T=1$ year. The sampling frequencies we consider range
from $\Delta_{n}=1$ second to $5$~minutes. The highest sampling frequencies
serve to validate the asymptotic theory contained above, while the lower
frequencies serve as proxies for situations where sparse sampling is employed
as a means to reduce the adverse impact of market microstructure
noise.\looseness=-1

The tables and graphs that follow report the results of $5\mbox{,}000$ simulations.
In order to validate the asymptotic theory developed above, we start
with the
highest sampling frequency of $\Delta_{n}=1$ second. Below, we will examine
the accuracy of the tests as a function of the sampling frequencies ranging
from $5$ seconds to $5$ minutes. The data generating process is the stochastic
volatility model $dX_{t}=\sigma_{t}\,dW_{t}+\theta dY_{t},$ with $\sigma
_{t}=v_{t}^{1/2}$, $dv_{t}=\chi(\eta-v_{t})\,dt+\xi v_{t}^{1/2}\,dB_{t}+dJ_{t}$,
$\mathbb{E}(W_{t}B_{t})=\overline{\rho} t$, $\eta^{1/2}=0.25,$ $\xi=0.5,$
$\chi=5$, $\overline{\rho}=-0.5,$ $J$ is a compound Poisson jump
process with
jumps that are uniformly distributed on $[-30\%,30\%]$ and $X_{0}=1$.
The jump
process $Y$ is either a $\beta$-stable process with $\beta=1$, that is, a
Cauchy process (which has infinite activity, and will be our model under
$\Omega_{T}^{i\beta}$) or a compound Poisson process (which has finite
activity, and will be our model under $\Omega_{T}^{f}$). In the
infinite-activity case, empirical estimates of $\beta$ reported in
\citet{yacjacod09b}
are higher than $1;$ simulation results using $\beta$-stable processes with
such values are qualitatively similar. In the finite-activity case, the jump
size of the compound Poisson process is drawn from a truncated Normal
distribution with mean $0$ and standard deviation $0.10$, designed to produce
jumps greater in magnitude than $0.05$. The estimator is implemented
with a
truncation rate $u_{n}=\alpha\Delta_{n}^{\varpi}$ where $\varpi=0.20$.
In the
results, we report the level of truncation $u_{n}$ indexed by the number
$\alpha$ of standard deviations of the continuous martingale part of the
process, defined in multiples of the long-term volatility parameter
$\eta^{1/2}.$ That is, $\alpha$ is such that $u_{n}=\alpha\eta
^{1/2}\Delta
_{n}^{1/2}.$

The statistic $S_{n}$ is implemented with $k=2$ and $p=4.$ When $p=4$, the
constants appearing in (\ref{T-9}) are%
%
\begin{equation}
N(4,k)=\frac{16}{3}k ( 2k^{2}-k-1 ) . \label{eqN4,k}%
\end{equation}
For the second test, $S_{n}^{\prime}$ is implemented with $p^{\prime}=4,$
$p=3$ and the second truncation level $\gamma u_{n}$ related to the first
according to $\alpha^{\prime}=2\alpha,$ that is, $\gamma=2.$

Given $\eta,$ the scale parameter $\theta$ [or equivalently $A_{t}$ in
(\ref{BA2}), which here is a constant] of the stable process in
simulations is
calibrated to deliver different various values of the tail probability
$\mathit{TP}=\mathbb{P}(|\Delta Y_{t}|\geq4\eta^{1/2}\Delta_{n}^{1/2})$ reported
in the
tables as low ($\mathit{TP}=0.01$), medium ($\mathit{TP}=0.05$) and high ($\mathit{TP}=0.10$). For the
Poisson process, it is the value of the arrival rate parameter $\lambda
$ that
is set to generate the desired level of jump tail probability as low
($\lambda=2/23\mbox{,}400$),\vadjust{\goodbreak} medium ($\lambda=10/23\mbox{,}400$) or high ($\lambda
=50/23\mbox{,}400$). In the various simulations' design, we hold $\eta$ fixed.
Therefore the tail probability parameter controls the relative scale of the
jump component of the semimartingale relative to its continuous counterpart.

\begin{table}
\def\arraystretch{0.9}
\tabcolsep=0pt
\caption{Monte Carlo rejection rate for the test of finite activity using the
statistic $S_n$ (testing~$ H_{0}\dvtx \Omega_{T}^{f}\cap\Omega_{T}^{W}$  vs.  $H_{1}\dvtx \Omega_{T}^{i\beta}$ using $S_{n}$)}
\label{tableMCS1H0}
\begin{tabular*}{\textwidth}{@{\extracolsep{4in minus 4in}}lcccccccc@{}}
\hline
\multirow{3}{34pt}[2.7pt]{\textbf{Finite jump intensity}} & \multirow{3}{41pt}[2.7pt]{\centering\textbf{Test theoretical level}} & \multicolumn{7}{c@{}}{\textbf{Test empirical level for a degree of truncation} $\bolds\alpha$}\\[-5pt]
  &   & \multicolumn{7}{c@{}}{\hrulefill}\\
  &   & $\bolds{6}$ & $\bolds{7}$ & $\bolds{8}$ & $\bolds{9}$ & $\bolds{10}$ & $\bolds{12}$ & $\bolds{15}$ \\
\hline
Low & $10\%$ & \hphantom{1}{$9.6\%$} & \hphantom{1}$9.6\%$ & \hphantom{1}$9.5\%$ & \hphantom{1}$9.5\%
$ &
\hphantom{1}$9.5\%$ & \hphantom{1}$9.5\%$ & \hphantom{1}$9.5\%$\\
& \hphantom{1}$5\%$ & \hphantom{1}$5.6\%$ & \hphantom{1}$4.5\%$ & \hphantom{1}$4.4\%$ & \hphantom{1}$4.4\%$ & \hphantom{1}$4.4\%$ & \hphantom{1}$4.4\%$ &
\hphantom{1}$4.4\%$\\[3pt]
Medium & $10\%$ & {$11.5\%$} & $10.7\%$ & $10.5\%$ &
$10.5\%$ & $10.5\%$ & $10.4\%$ & $10.4\%$\\
& \hphantom{1}$5\%$ & \hphantom{1}$6.1\%$ & \hphantom{1}$5.0\%$ & \hphantom{1}$5.0\%$ & \hphantom{1}$5.0\%$ & \hphantom{1}$5.0\%$ & \hphantom{1}$5.0\%$ &
\hphantom{1}$5.0\%$\\[3pt]
High & $10\%$ & {$11.8\%$} & $10.4\%$ & $10.3\%$ &
$10.3\%$
& $10.3\%$ & $10.2\%$ & $10.2\%$\\
& \hphantom{1}$5\%$ & \hphantom{1}$5.9\%$ & \hphantom{1}$5.1\%$ & \hphantom{1}$5.0\%$ & \hphantom{1}$5.0\%$ & \hphantom{1}$5.0\%$ & \hphantom{1}$5.0\%$ &
\hphantom{1}$5.0\%$\\
\hline
\end{tabular*}  \vspace*{-3pt}
\end{table}
%

\begin{figure}[b]
\vspace*{-3pt}
\includegraphics{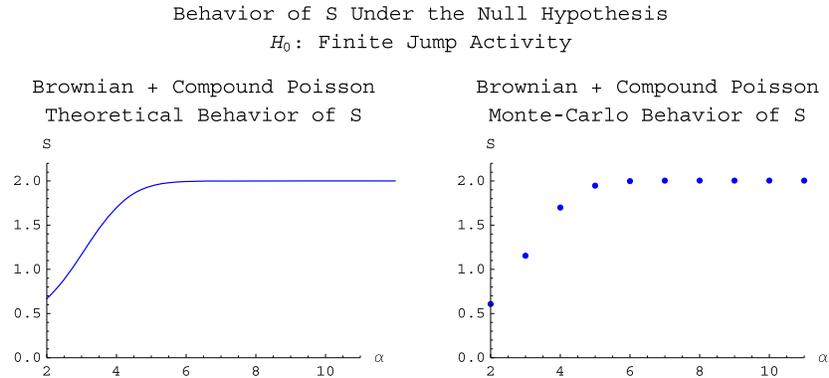}

\caption{Behavior of the test statistic $S_{n}$ under $H_{0}\dvtx \Omega_{T}%
^{f}\cap\Omega_{T}^{W}$ for varying degrees of truncation $\alpha$.}%
\label{figS1H0Curves}
\end{figure}

Table~\ref{tableMCS1H0} reports the Monte Carlo rejection rates of the test
of\vspace*{-1pt} $H_{0}\dvtx \Omega_{T}^{f}\cap\Omega_{T}^{W}$ vs. $H_{1}\dvtx \Omega
_{T}^{i\beta}$ at
the $10\%$ and $5\%$ level, using the test statistic $S_{n},$ for various
levels of truncation $\alpha$. Recall that for concreteness $\alpha$ is
expressed as a number of standard deviations of the Brownian part of
$X$. We
find that the test behaves well, with empirical test levels close to their
theoretical counterparts, and so for a wide range of values of $\alpha.$

The limit in probability of $S_{n}$ is $k^{p/2-1}=2$ with our choice of $k$
and $p$, provided $u_{n}/\sqrt{\Delta_{n}}\rightarrow\infty$. With our
notation $\alpha$ denoting the threshold expressed as the number of
(normalized) standard deviations of the Brownian part, this means that this
limit holds when $\alpha=\alpha_{n}$ goes to infinity. In practice we choose
$\alpha$ between $6$ and $15$, and this introduces a bias. To evaluate this
bias, in Figure~\ref{figS1H0Curves}, we\vadjust{\goodbreak} plot as a function of $\alpha$ the
limiting value of $S_{n}$ under $H_{0},$ when $\alpha$ stays fixed
instead of
increasing to infinity, and compare it with the theoretical limit $2$
(in the
simple model used for the simulations it is possible to numerically
compute it
with any desired accuracy); this is the left graph, and on the right
graph we
draw the corresponding average value of $S_{n}$ from the Monte Carlo
simulations, both as functions of $\alpha$. It is worth noting that the
behavior of $S_{n}$ is driven by that of the Brownian motion component, since
the truncation effectively eliminates the finite-activity jumps. For very
small values of $\alpha$, $\alpha<5$, $S_{n}$ starts below $2$ as
predicted by
the theoretical behavior of the continuous martingale part of $X.$ Once
$\alpha$ gets above $5,$ the theoretical limit of $S_{n}$ will then
remain at
$2$ as long as $\alpha$ is not so large as to start including some jumps.

\begin{figure}

\includegraphics{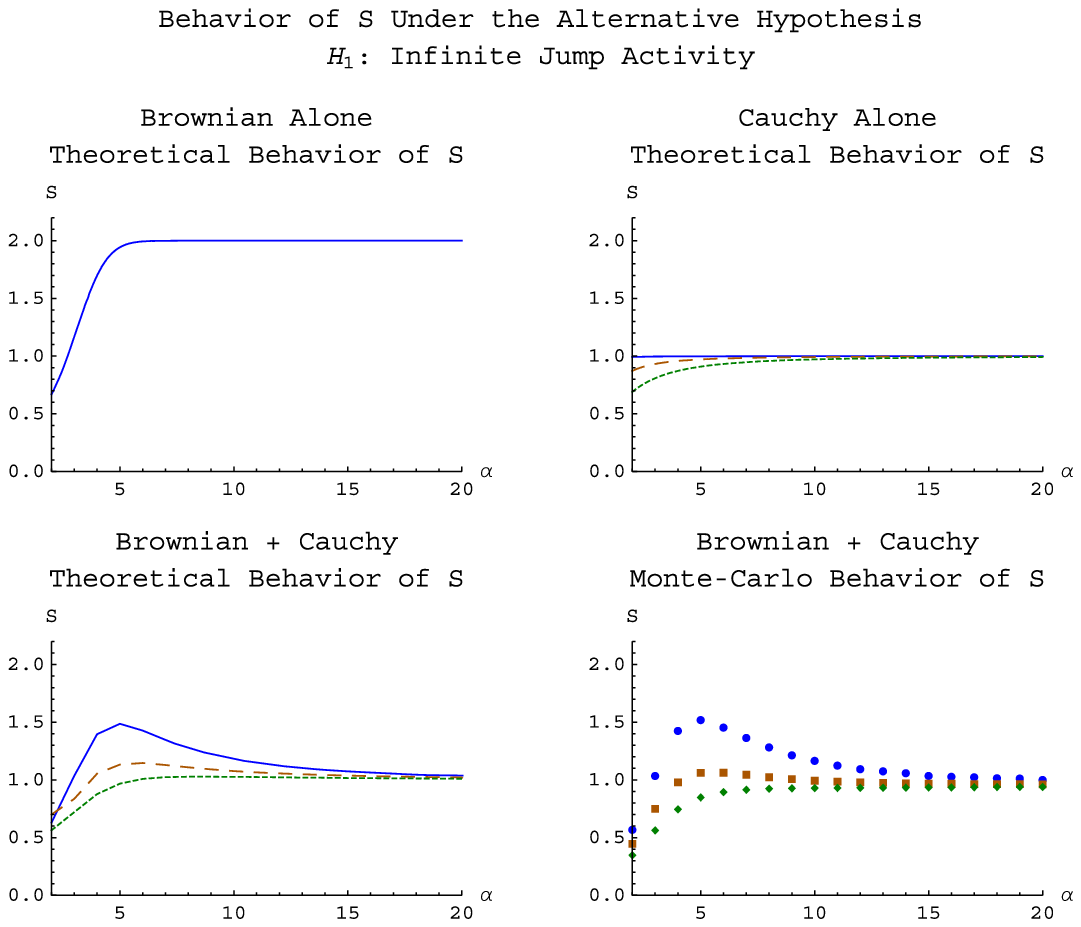}

\caption{Behavior of the test statistic $S_{n}$ under $H_{1}\dvtx \Omega
_{T}^{i\beta}$ for varying degrees of truncation~$\alpha$.}
\label{figS1H1Curves}
\end{figure}

Figure~\ref{figS1H1Curves} plots the limiting value of $S_{n}$ (as a function
of $\alpha$ again) under~$H_{1},$ where the model consists of a
Brownian and a
Cauchy components as described above. The upper left graph is the theoretical
behavior of $S_{n}$ for the Brownian component taken in isolation, the upper
right graph the theoretical behavior for the Cauchy component taken in
isolation, the lower left graph the theoretical behavior for the sum of the
two components, that is, the model we simulate from, and the lower
right graph
the corresponding results from the Monte Carlo simulations. When the Cauchy
process is present, the three curves in the figure correspond to the low
(solid curve), medium (long dashes) and high (short dashes) scale of
the jump
component relative to the continuous component. For very small values of
$\alpha,$ such as $\alpha<5,$ the behavior of $S_{n}$ tracks that of the
Brownian component, which is increasing in $\alpha$ toward $2.$ As
$\alpha$
increases, the curves then start reversing course and tend to $1$
instead, the
limit driven by the infinite-activity jump process as shown in (\ref{T-7}).
The higher the scale of the jump component relative to the continuous
component, the more the curve in the lower left graph approximates the
corresponding one in the upper right graph, that is, the more it resembles
that of a pure Cauchy process. For the Monte Carlo results, the low,
medium and high jump scales are represented by circles, squares and diamonds,
respectively. In all three cases, they track the predicted theoretical limits
closely for each value of the truncation level~$\alpha.$%

\begin{figure}

\includegraphics{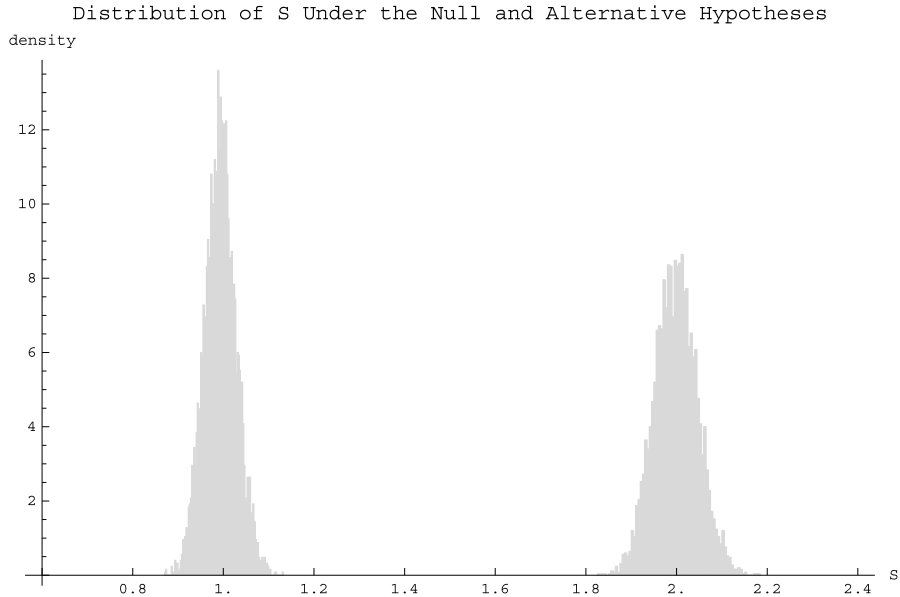}

\caption{Monte Carlo distributions of the unstandardized test statistic
$S_{n}$ under\vspace*{-1pt} $H_{0}\dvtx \Omega_{T}^{f}\cap\Omega_{T}^{W}$ (right
histogram) and
$H_{1}\dvtx \Omega_{T}^{i\beta}$ (left histogram).}%
\label{figS1H0H1HistNonStd}
\end{figure}

\begin{figure}

\includegraphics{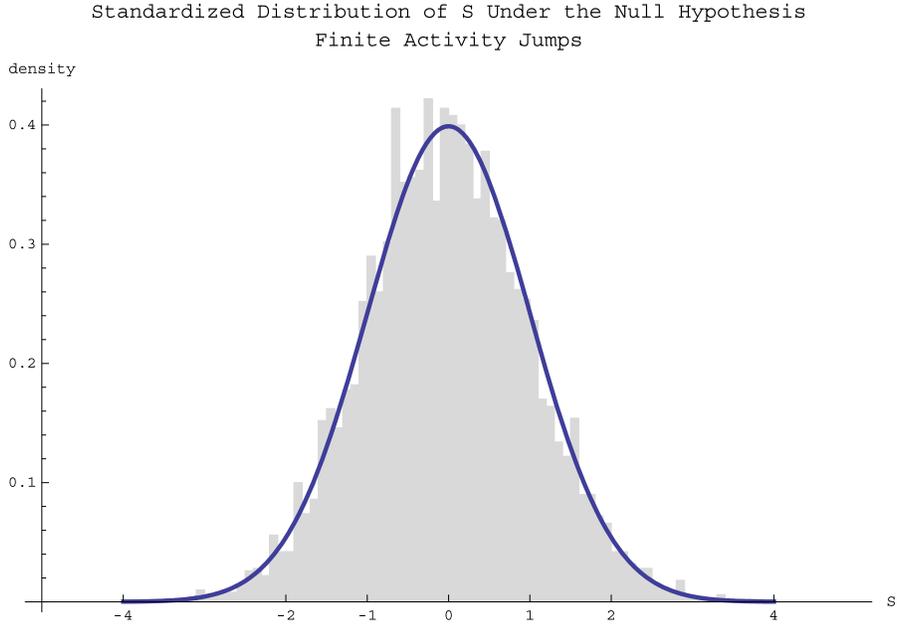}

\caption{Monte Carlo and asymptotic distributions of the standardized test
statistic $S_{n}$ under $H_{0}\dvtx \Omega_{T}^{f}\cap\Omega_{T}^{W}$.}%
\label{figS1H0HistStd}
\end{figure}

\begin{table}[b]
\tabcolsep=0pt
\tablewidth=260pt
\caption{Monte Carlo rejection rate for the test of infinite jump activity
using the statistic $S_{n}%
^{\prime}$ (testing~$ H_{0}\dvtx \Omega_{T}^{i\beta}\mbox{ vs. }H_{1}\dvtx \Omega_{T}^{f}$
using $S_{n}^{\prime}$)}
\label{tableMCS2H0}
\begin{tabular*}{260pt}{@{\extracolsep{\fill}}lcc@{}}
\hline
\textbf{Infinite}&\textbf{Test}&\textbf{Test}\\
\textbf{jump activity}&\textbf{theoretical level}&\textbf{empirical level}\\
\hline
Low & $10\%$ & \hphantom{1}$8.2\%$\\
& \hphantom{1}$5\%$ & \hphantom{1}$3.6\%$\\[3pt]
Medium & $10\%$ & \hphantom{1}$8.9\%$\\
& \hphantom{1}$5\%$ & \hphantom{1}$3.8\%$\\[3pt]
High & $10\%$ & $10.0\%$\\
& \hphantom{1}$5\%$ & \hphantom{1}$4.9\%$\\
\hline
\end{tabular*}
\end{table}

\begin{figure}

\includegraphics{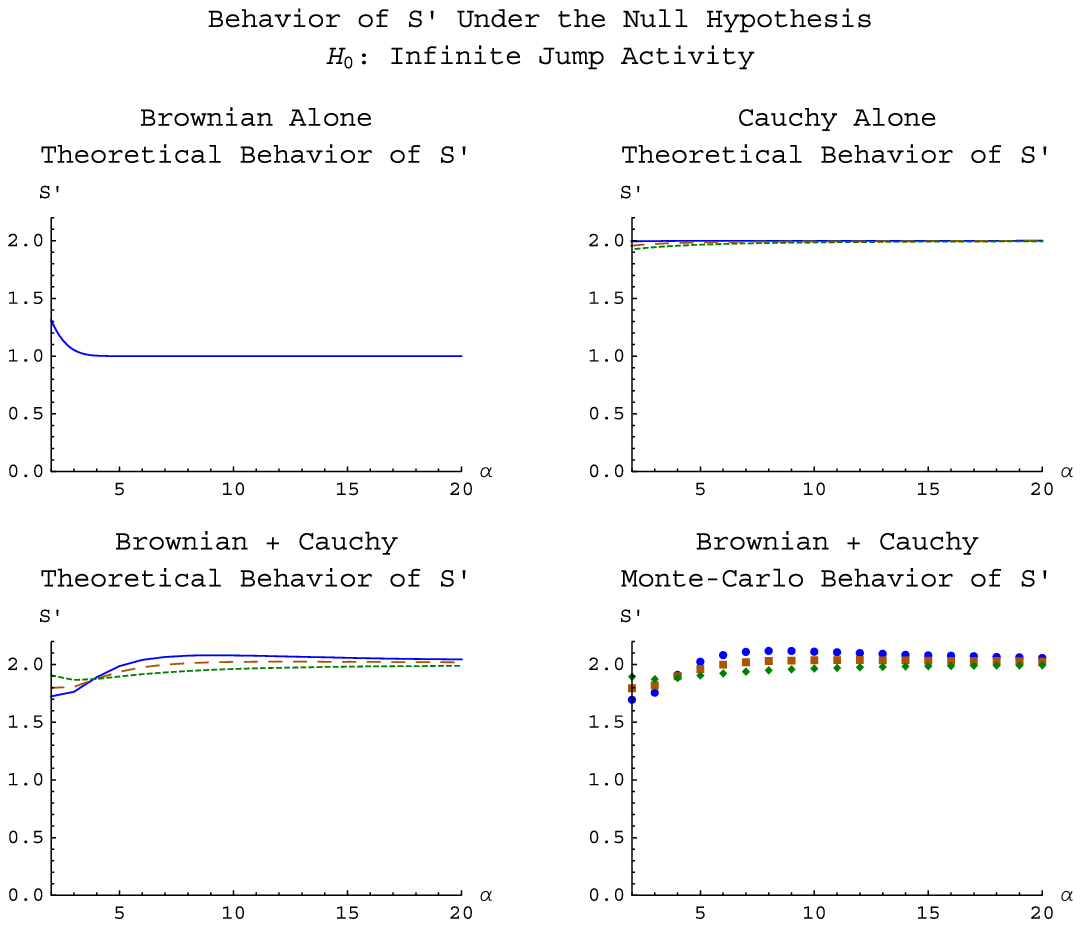}

\caption{Behavior of the test statistic $S_{n}^{\prime}$ under $H_{0}\dvtx \Omega_{T}^{i\beta}$ for varying degrees of truncation~$\alpha$.}%
\label{figS2H0Curves}%
\end{figure}

Finally, we report in Figure~\ref{figS1H0H1HistNonStd} histograms of the
values of the unstandardized $S_{n}$ computed under $H_{0}$, centered
in $2$
as expected from (\ref{T-6}), and $H_{1},$ centered in $1$ as expected from
(\ref{T-7}), respectively. Figure~\ref{figS1H0HistStd} reports the Monte
Carlo distribution of the statistic $S_{n},$ standardized according to
(\ref{eqCLTS36}), compared to the limiting $\mathcal{N}(0,1) $ distribution.

Next, we turn to the symmetric problem, that of testing of $H_{0}\dvtx \Omega
_{T}^{i\beta}$ vs. $H_{1}\dvtx \Omega_{T}^{f}\cap\Omega_{T}^{W},$ which we
do using
the statistic $S_{n}^{\prime}.$ The results are reported in Table
\ref{tableMCS2H0} for the test rejection rate under the null hypothesis.
Figure~\ref{figS2H0Curves} shows the limiting value of $S_{n}^{\prime}$ under
$H_{0},$ as expected from the theoretical limit of $\gamma^{p^{\prime}-p}=2$
given in (\ref{T-12}) (lower left graph for the Brownian plus Cauchy model)
and the corresponding average value of $S_{n}^{\prime}$ from the Monte Carlo
simulations (lower right graph). The upper graphs in the figure
correspond to
the Brownian alone (left, with a limit of~$1$) and Cauchy alone (right,
with a
limit of $2$) situations. As was the case for the previous test, higher values
of the jump scale parameter make the resulting semimartingale model
approximate the behavior of its infinite-activity jump component more
closely, while
for lower scale parameters and low values of $\alpha,$ the behavior is
determined by the continuous component. In all cases, we find that the Monte
Carlo results track the predicted theoretical limits for each value of the
truncation level $\alpha.$ Compared to Figure~\ref{figS1H1Curves}, we find
that the limit is approached with less precision than was the case of $S_{n},$
and requires larger values of $\alpha$ than was the case for the other test
when the infinite-activity jump process is mixed with the continuous
component.%

%

\begin{figure}

\includegraphics{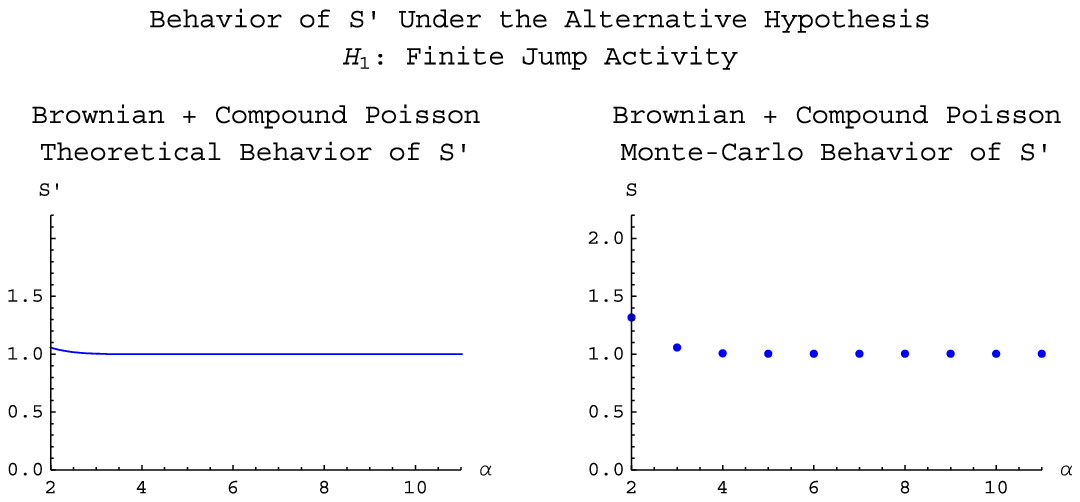}

\caption{Behavior of the test statistic $S_{n}^{\prime}$ under $H_{1}\dvtx \Omega_{T}^{f}\cap\Omega_{T}^{W}$ for varying degrees of truncation
$\alpha$.}%
\label{figS2H1Curves}
\end{figure}

\begin{figure}[b]

\includegraphics{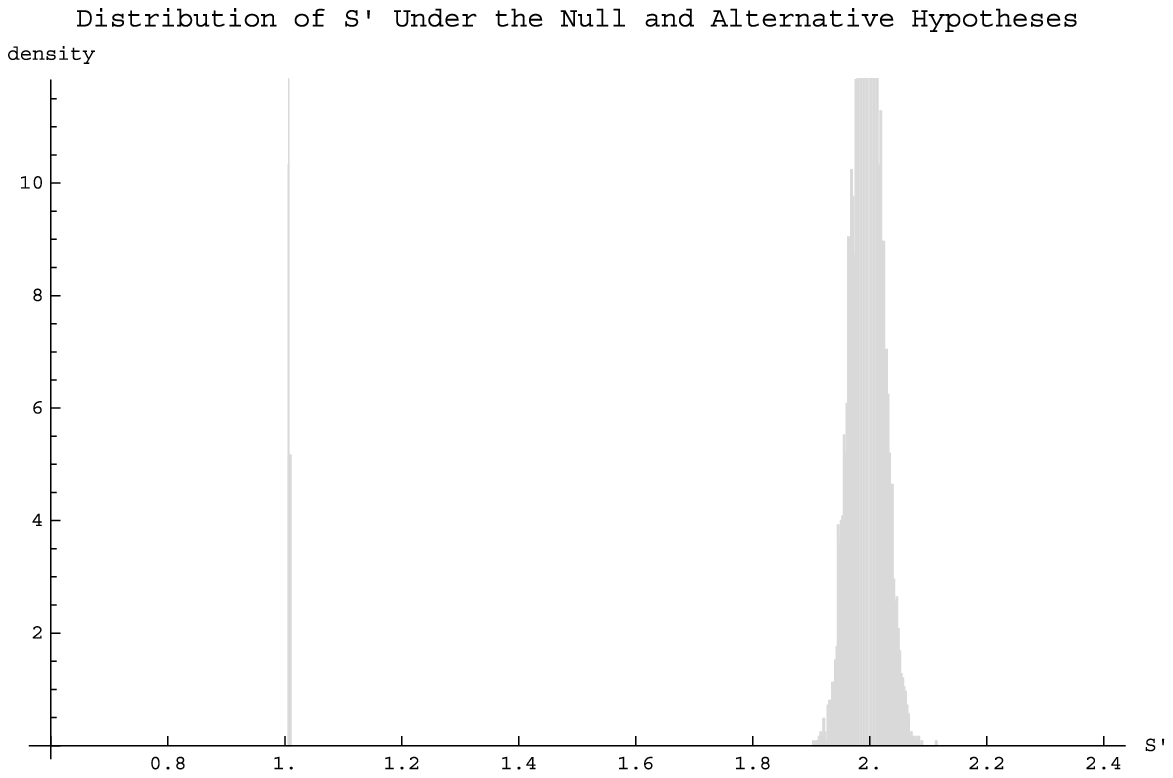}

\caption{Monte Carlo distributions of the unstandardized test statistic
$S_{n}^{\prime}$ under $H_{0}\dvtx \Omega_{T}^{i\beta}$\vspace*{-1pt} (right histogram) and
$H_{1}\dvtx \Omega_{T}^{f}$ (left histogram).}%
\label{figS2H0H1HistNonStd}
\end{figure}

Figure~\ref{figS2H1Curves} shows the point limit results for
$S_{n}^{\prime}$
under $H_{1}.$ There, the theoretical limit of $1$ from (\ref{T-13}) is
reached quickly, as was the case in Figure~\ref{figS1H0Curves}. This is the
case provided that the values of $\alpha$ are not so large that jumps start
being retained in the calculations. In particular, the value of $S_{n}%
^{\prime}$ can become large if $\alpha$ is such that a jump is just included
at the $u_{n}$ cutoff but not yet at the $\gamma u_{n}$ cutoff. But in most
cases the statistic is very close to $1$ when the number of jumps is finite,
simply because once $n$ is large enough for the truncation level to have
eliminated all the jumps at the higher level of truncation $\gamma u_{n},$
then there are very few Brownian increments between truncation levels $u_{n}$
and $\gamma u_{n}$; therefore $B(p,\gamma u_{n},\Delta_{n})_{T}$ is only
marginally larger than $B(p,u_{n},\Delta_{n})_{T},$ and similarly for power
$p^{\prime}$.

\begin{figure}

\includegraphics{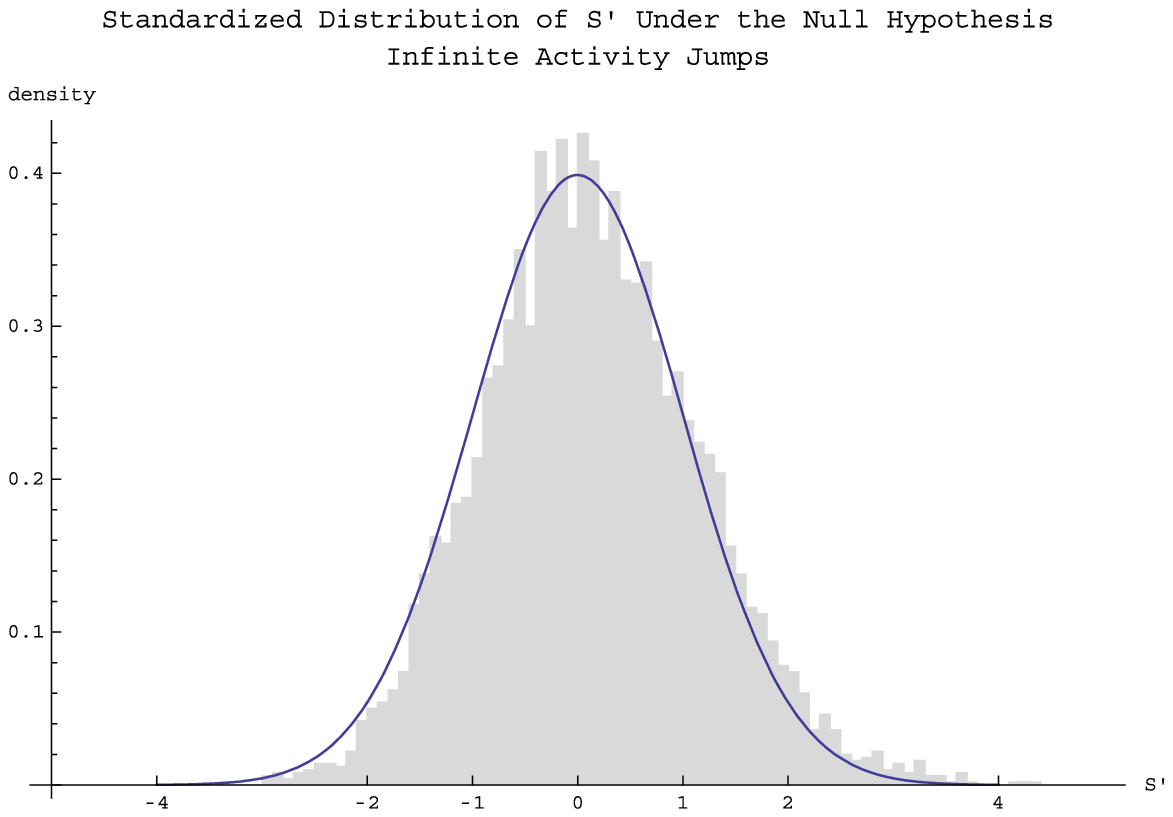}

\caption{Monte Carlo and asymptotic distributions of the standardized test
statistic $S_{n}^{\prime}$ under $H_{0}\dvtx \Omega_{T}^{i\beta}$.}%
\label{figS2H0HistStd}
\end{figure}

Figure~\ref{figS2H0H1HistNonStd} reports the Monte Carlo distributions of
$S_{n}^{\prime}$ under $H_{0}$ and $H_{1};$ they are centered at $2$
and $1$
as expected. Under $H_{1},$ we note that $S_{n}^{\prime}$ displays very little
variability with the provision above regarding the value of $\alpha.$
As a
general rule, the test based on $S_{n}^{\prime}$ appears more sensitive to
values of $\alpha$ than the test based on $S_{n}.$ Figure
\ref{figS2H0HistStd} reports the Monte Carlo and asymptotic
distribution of
the standardized $S_{n}^{\prime}$ under $H_{0},$ according to (\ref{eqCLTS43}).

\begin{figure}

\includegraphics{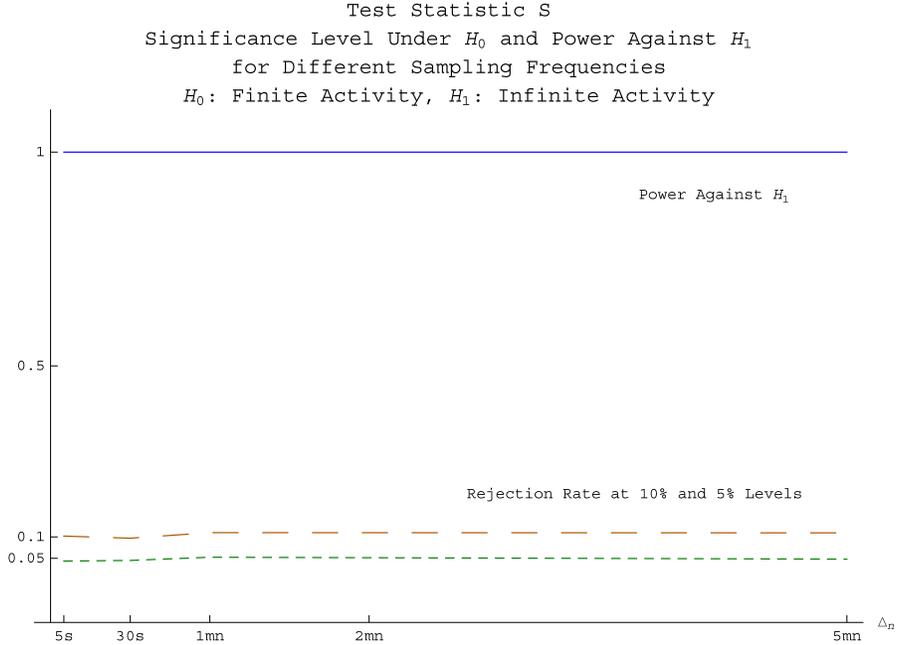}

\caption{Level and power of the test based on the statistic $S$ as a function
of the sampling frequency.}%
\label{figS1AllDel}
\end{figure}

\begin{figure}

\includegraphics{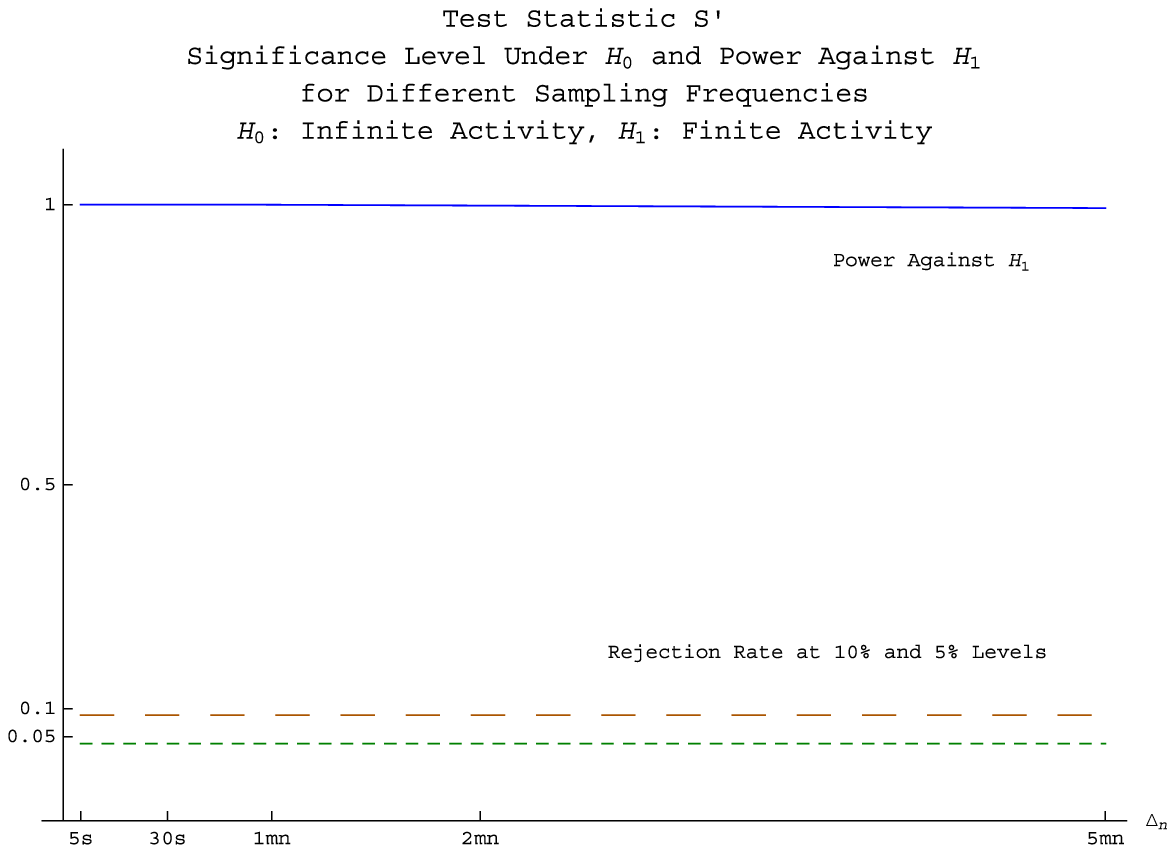}

\caption{Level and power of the test based on the statistic $S^{\prime}$ as a
function of the sampling frequency.}%
\label{figS2AllDel}
\end{figure}

Finally, we examine the accuracy of the tests as a function of the sampling
frequencies. We consider for this purpose two experiments. In the first one,
we fix the observation length to $T=5$ days (one week) and consider sampling
frequencies from $5$ seconds to $5$ minutes. The sample size $n$ therefore
decreases by a factor $60$ over the range of values of $\Delta_{n}$
considered. In the second experiment, we consider the same sampling
frequencies but increase the length of the observation window $T$ from $5$
days to $1$ year, in order to keep the sample size approximately
constant over
the range of values of $\Delta_{n}.$ The results of the first
experiment are
reported in Figures~\ref{figS1AllDel} and~\ref{figS2AllDel},
respectively. We
find little size distortion or power loss. The results of the second
experiment are, not surprisingly, better and are not shown here in
order to
save space.

\section{Empirical results}\label{secdata}%

In this section, we apply our two test statistics to real data,
consisting of
all transactions recorded during the year 2006 on two of the most actively
traded stocks, Intel (INTC) and Microsoft (MSFT). The data source is
the TAQ
database. Using the correction variables in the data set, we retain only
transactions that are labeled ``good
trades''
by the exchanges: regular trades that were not corrected, changed, or
signified as canceled or in error; and original trades which were later
corrected, in which case the trade record contains the corrected data
for the
trade. Beyond that, no further adjustment to the data is made.

\begin{figure}

\includegraphics{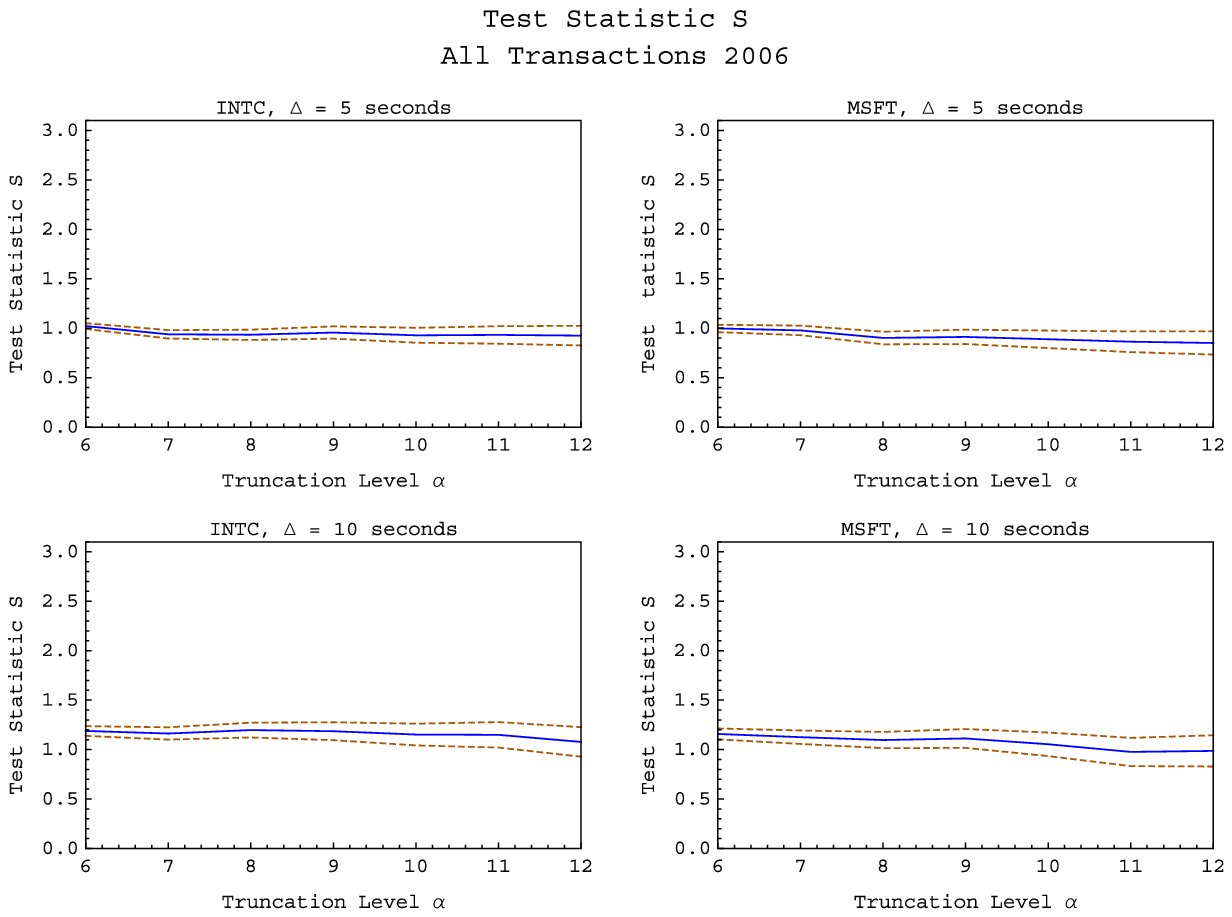}

\caption{Values of the test statistic $S_{n}$ computed for all 2006
transactions of INTC and MSFT, sampled at $5$ and $10$ second
intervals.}%
\label{figINTCMSFTS1}
\end{figure}

\begin{figure}

\includegraphics{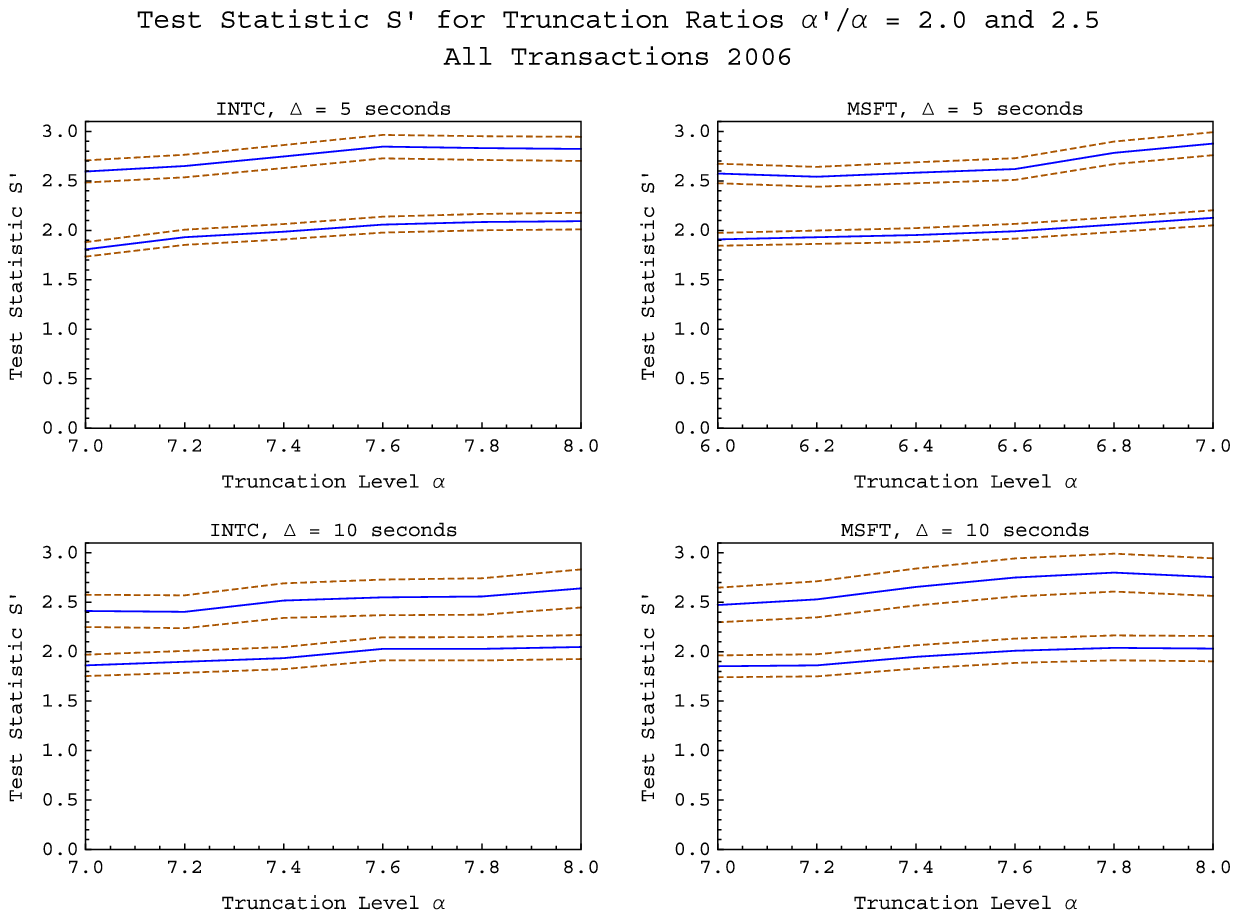}

\caption{Values of the test statistic $S_{n}^{\prime}$ computed for all 2006
transactions of INTC and MSFT, sampled at $5$ and $10$ second
intervals.}%
\label{figINTCMSFTS2}
\end{figure}

We first consider the test where the null hypothesis consists of finite jump
activity. Figure~\ref{figINTCMSFTS1} shows the values of the test statistic
$S_{n},$ along with a $95\%$ confidence interval computed from the asymptotic
distribution (\ref{eqCLTS36}) under the finite-activity null hypothesis,
plotted for a range of values of the truncation index $\alpha$. The truncation
level plays the same role as that of a bandwidth parameter in classical
nonparametric estimator, and it is therefore important that it be properly
scaled initially. For this purpose, the values of $\alpha$ are indexed in
terms of standard deviations of the continuous martingale part of the
log-price: we first estimate the volatility of the continuous part of $X$
using the sum of squared increments that are smaller than $\Delta_{n}^{1/2}$
(meaning that we would retain the increments of up to four standard deviations
if the annualized volatility of the stock were $25\%$ per year) and
then use
that estimate to form the initial cutoff level used in the construction
of the
test statistic. To account for potential time series variation in the
volatility process $\sigma_{t}$, that procedure is implemented
separately for
each day and we compute the sum, for that day, of the increments that are
smaller than the cutoff, to the appropriate power $p$ required by the test
statistic. For the full year, we then add the power variations computed for
each day. We then compute the results corresponding to a range of
values of
$u_{n}$ indexed by $\alpha$.

In order to account for the presence of market microstructure noise in the
data, we compute the limiting values of our test statistics in the case where
the noise is a pure additive noise; the limit is then $1/k$,
independent of
$p$. We also consider the dependence of the test statistics as the sampling
interval increases, and the signal-to-noise ratio presumably improves.

The test statistic is implemented with $p=4$ and $k=2,$ using $\Delta_{n}=5$
seconds in the upper panels, and $\Delta_{n}=10$ seconds in the lower panels.
As a result, from Theorem~\ref{TT1}, $S_{n}$ should go to $k^{p/2-1}=2$ under
the null of finite activity, and to $1$ under the alternative of infinite
activity. As the plots show, we find that $S_{n}$ is close to $1,$
which leads
us to reject the null hypothesis of finite activity.

Next, we turn to the test where the null hypothesis consists of
infinite jump
activity. Figure~\ref{figINTCMSFTS2} shows the values of the test statistic
$S_{n}^{\prime},$ along with a $95\%$ confidence interval computed from the
asymptotic distribution (\ref{eqCLTS43}) under the infinite-activity null
hypothesis, plotted for a range of values of the truncation index
$\alpha$.
The two curves correspond to values of $\gamma=\alpha^{\prime}/\alpha$ equal
to $2.0$ and $2.5,$ respectively. The test statistic is implemented with
$p^{\prime}=4$ and $p=3,$ using $\Delta_{n}=5$ seconds in the upper panels,
and $\Delta_{n}=10$ seconds in the lower panels. Recall from Theorem~\ref{TT3}
that under the null of infinite activity, $S_{n}^{\prime}$ should go to
$\gamma^{p^{\prime}-p}=2.0$ and $2.5$, respectively, for the two
curves, and to
$1$ under the alternative of finite activity. We find that
$S_{n}^{\prime}$ is
close to the predicted value $\gamma^{p^{\prime}-p},$ which leads us to not
reject the null hypothesis of infinite activity.

\begin{figure}

\includegraphics{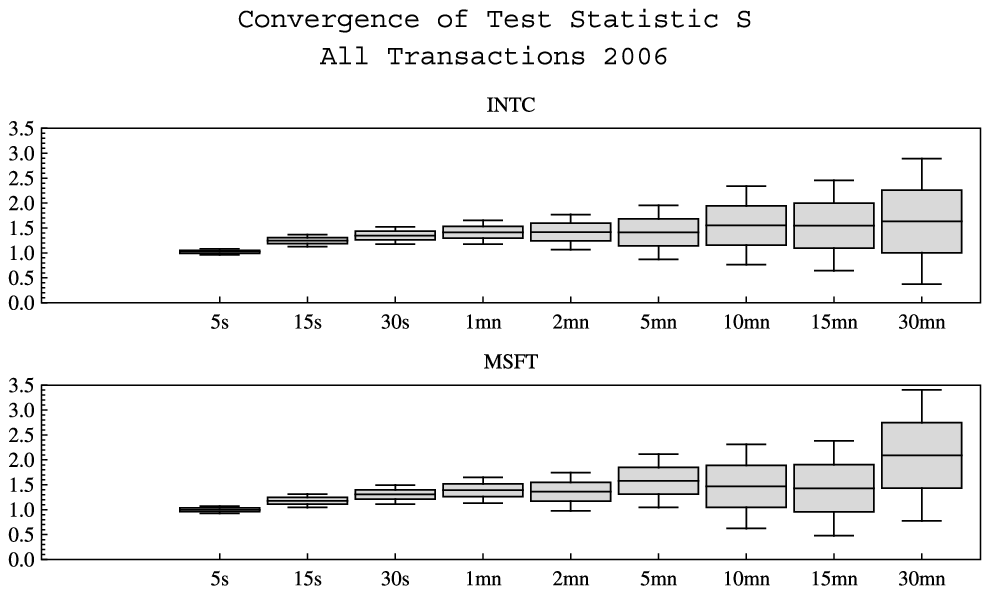}

\caption{Box plot for the test statistic $S_{n}$ computed for all 2006
transactions of INTC and MSFT, sampled at time intervals ranging from $5$
seconds to $30$ minutes.}%
\label{figINTCMSFTS1BoxPlotConvergence}
\end{figure}

\begin{figure}

\includegraphics{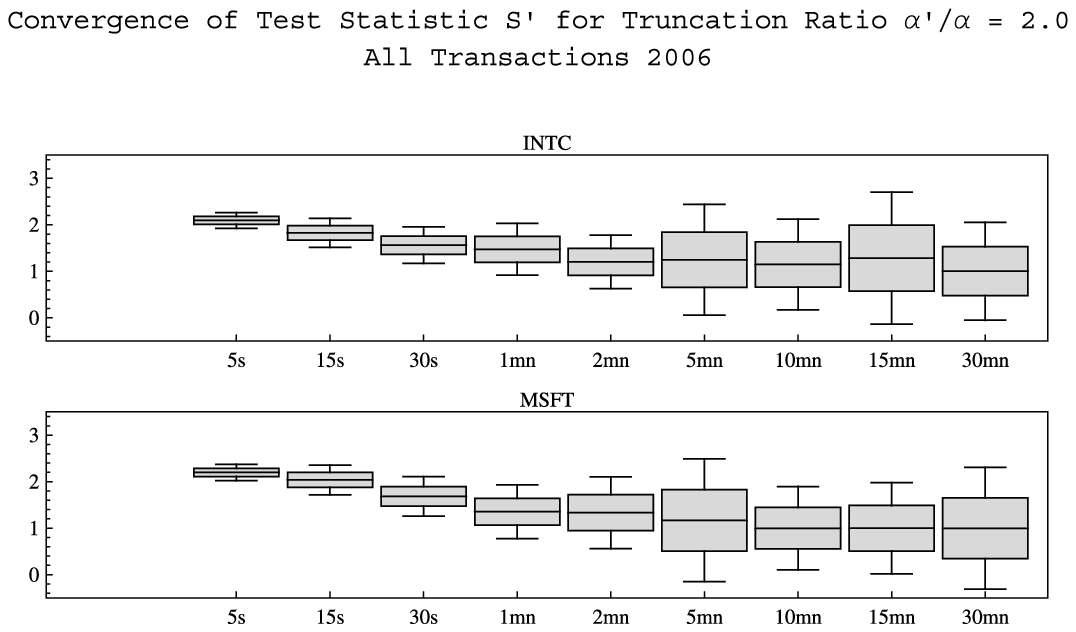}

\caption{Box plot for the test statistic $S_{n}^{\prime}$ computed for all
2006 transactions of INTC and MSFT, sampled at time intervals ranging
from $5$
seconds to $30$ minutes.}%
\label{figINTCMSFTS2BoxPlotConvergence}
\end{figure}

To summarize, the answer from both tests appears indicative of infinite jump
activity in those data: using $S_{n},$ we reject the null of finite activity,
while using $S_{n}^{\prime}$ we do not reject the null of infinite activity.
Finally, we illustrate in Figure~\ref{figINTCMSFTS1BoxPlotConvergence} the
convergence of $S_{n}$ to $1$ as the sampling interval decreases, indicating
that the null hypothesis of finite activity is rejected when high-frequency
data (of the order of seconds) are used. On the other hand, we see that
using longer sampling intervals (of the order of minutes) makes it
impossible to
reject the null hypothesis of finite activity using $S_{n}.$ This is
compatible with the fact that small jumps occurring over short time intervals
can be aggregated or smoothed out over longer time intervals. For the second
testing situation, Figure~\ref{figINTCMSFTS2BoxPlotConvergence} shows the
convergence of $S_{n}^{\prime}$ to $2$ as the sampling interval
decreases from
$30$ minutes to $5$ seconds, indicating that the null hypothesis of infinite
activity is not rejected at high frequency. As in the first test, lower
frequency data tend to be more compatible with finite jump activity. In both
cases, longer sampling intervals $\Delta_{n}$ over the same sampling length
$T$ lead to a reduction in the sample size $n=T/\Delta_{n},$ which generally
leads to an increase in the variance of the test statistic, making it more
difficult ceteris paribus to reject the null hypothesis.

\section*{Acknowledgments}
We are very grateful to a referee, an Associate Editor and the Editor
for many
helpful comments.

\begin{supplement}
\stitle{Supplement to ``Testing whether jumps have finite or infinite activity''}
\slink[doi]{10.1214/11-AOS873SUPP} 
\sdatatype{.pdf}
\sfilename{aos873\_suppl.pdf}
\sdescription{This supplementary article contains a few additional technical details about the assumptions made in this paper, and the proofs of all results.}
\end{supplement}

%

\printaddresses

\end{document}